\documentclass[11pt]{article}
\usepackage{hyperref}
\usepackage{graphicx}
\usepackage[leqno]{amsmath}
\usepackage{amssymb,amsthm,mathabx}	
\usepackage{array}
\usepackage{float}
\restylefloat{figure,wrapfig}
\usepackage[numbers]{natbib}
\newskip\stdskip                      
\hypersetup{%
  bookmarksnumbered=true,%
  bookmarks=true,%
  colorlinks=true,%
  linkcolor=blue,%
  citecolor=blue,%
  filecolor=blue,%
  menucolor=blue,%
  pagecolor=blue,%
  urlcolor=blue,%
  pdfnewwindow=true,%
  pdfstartview=FitBH}
\let\fullref\autoref
\newtheorem*{theorem*}{Theorem}
\newtheorem*{corollary*}{Corollary}
\newtheorem{theorem}{Theorem}[section]
\newtheorem{proposition}[theorem]{Proposition}
\newtheorem{corollary}[theorem]{Corollary}
\newtheorem{lemma}[theorem]{Lemma}
\newtheorem{definition}[theorem]{Definition}
\newtheorem{example}[theorem]{Example}
\newtheorem{remark}[theorem]{Remark}
\numberwithin{equation}{section}
\numberwithin{figure}{section}
\def\makeautorefname#1#2{\expandafter\def\csname#1autorefname\endcsname{#2}}
\makeautorefname{equation}{Equation}%
\makeautorefname{footnote}{footnote}%
\makeautorefname{figure}{Figure}%
\makeautorefname{chapter}{Chapter}%
\makeautorefname{section}{Section}%
\makeautorefname{subsection}{Section}%
\makeautorefname{subsubsection}{Section}%
\def \cA{\mathcal A}  \def \cB{\mathcal B}  \def \cC{\mathcal C}  
\def \cH{\mathcal H}  \def \cL{\mathcal L}  
\def \cM{\mathcal M}  \def \cV{\mathcal V}	
\def \R{\ensuremath{\,\mathbb{R}}} 		\def \Z{\ensuremath{\,\mathbb{Z}}}  
\def \S{\ensuremath{\,\mathbb  {S}}} 	
\def \C{\ensuremath{\, \mathbb{C}}} 	
\def \Q{\ensuremath{\, \mathbb{Q}}} 	
\def \D{\hbox {\mathversion{bold}$D$\mathversion {normal}}}
\def \thm #1{		{\begin{theorem} 		{#1}\end{theorem}}\par		}
\def \prop #1{	{\begin{proposition}{#1}\end{proposition}}\par }
\def \cor #1{		{\begin{corollary} 	{#1}\end{corollary}}\par	}
\def \lem #1{		{\begin{lemma} 			{#1}\end{lemma}}\par		}
\def \defin #1{	{\begin{definition} {#1}\end{definition}}\par	}

\def \rem #1{		{\begin{remark} 		{#1}\end{remark}}\par			}
\def \Comp #1{[\hspace{-2pt}[#1]\hspace{-2pt}]}
\def \Ker{\mathop{\rm Ker}}
\def \Im{\mathop{\rm Im}}
\def \Fix{{\mathop{\rm Fix}}}
\def \Oplus{\mathop{\text{\Large $\oplus$}}}
\def \Otimes{\mathop{\text{\Large $\otimes$}}}
\def \cl{{\mathop{c\hskip 0.5pt l}}}
\def \op{{\mathop{op}}}
\def \st{{\mathop{St}}}
\def \Ga{{\mathit{\Gamma }}}
\let \eps=\varepsilon
 
\def \id{\text{\rm id}}
\def \cone{\mathop{\rm Cone}}
\def \downsquigarrow{\rotatebox[origin=c]{270}{$\rightsquigarrow$}}
\def \wtilde #1{\hspace{1pt} \widetilde {\hspace{-1pt}  #1}{}}
\def \wbar #1{\hspace{1pt} \overline {\hspace{-1pt}  #1}{}}

\def \noespace{\itemsep=0pt\topsep=0pt\partopsep=0pt\parsep=0pt\parskip=1pt}
\def \rightarrowvar #1{\hbox to #1 {\rightarrowfill}} 
\def \Ajvert #1#2{\vrule height #1pt depth #2pt width 0pt}
\let \dis=\displaystyle
\let \scr=\scriptstyle
\newlength\savedwidth

\begin{document}
\title{Khovanov homology for signed divides.}
\date{}
\author{Olivier Couture\\
\\
\small Institut de Math\'ematiques de Bourgogne\\
Universit\'e de Bourgogne\\
UFR Sciences et Techniques\\
9 avenue Alain Savary - B.P 47870\\
21078 DIJON Cedex - France\\
olivier.couture@u-bourgogne.fr}
\maketitle
\begin{flushright}
\emph{Dedicated to Bernard Perron}
\end{flushright}

\begin{abstract}
The purpose of this paper is to interpret polynomial invariants of strongly invertible links in 
terms of Khovanov homology theory. To a divide, that is a proper generic immersion of a finite number of copies of 
the unit interval and circles in a 2-disc, one can associate a strongly invertible link in the 3-sphere. 
This can be generalized to signed divides : divides with + or $-$ sign assignment to 
each crossing point. Conversely, to any link $L$ that is strongly invertible for an involution $j$, 
one can associate a signed divide. Two strongly invertible links that are isotopic 
through an isotopy respecting the involution are called strongly equivalent. Such isotopies give rise to
moves on divides. In a previous paper of the author \cite{C}, one can find an exhaustive list of moves that 
preserves strong equivalence, together with a polynomial invariant for these moves, giving therefore 
an invariant for strong equivalence of the associated strongly invertible links. We prove in this paper that 
this polynomial can be seen as the graded Euler characteristic of a graded complex of vector spaces. 
Homology of such complexes is invariant for the moves on divides and so is invariant through strong equivalence 
of strongly invertible links.
\end{abstract}

57M27

\section*{Introduction.}\label{intro} %
%
A \emph{divide} $\Ga$ is the image of a proper generic immersion of a finite number of intervals and circles into 
the unit 2-disc $\D^2$ of $\R^2$. With a divide $\Ga$, N.A'Campo \cite {AC} associates a link $\cL(\Ga)$ in the 
unit 3-sphere of the tangent space $T\D^2$:
$$
		\cL(\Ga)=\{ (p,v)\in T\D^2\ :\ p\in \Ga,\ v\in T_p\Ga,\ \|p\|^2+\|v\|^2=1\}.	
$$
This link has natural orientation and is strongly invertible with respect to the involution $j(p,v)=(p,-v)$.

In Couture--Perron \cite{CP} is given a generalization of divides. Let $(x,y)$ be coordinates in $\D^2$ such that 
the restriction to $\Ga$ of the projection $\pi_1:(x,y)\mapsto x$ is a Morse function. A \emph{Morse signed divide}
(MS--divide) relative to $\pi_1$ stands for such a divide  with $+$ or $-$ sign assignment to each double point of 
$\Ga$. Furthermore, if there exists $a\in ]0,1[$ such that all maxima (resp. minima) of ${\pi_1}_{\mid \Ga}$ 
project on $a$ (resp. $-a$) and all double points in $]-a;a[$, the MS--divide is called \emph{ordered (OMS--divide)}.

We also define a link associated with a MS--divide (see \cite{CP}), which is strongly 
invertible with respect to the involution $j(p,v)=(p,-v)$. If all signs are positive, this link is no more than 
the link of the divide without signs. The interest of OMS--divides rather than MS--divides is to obtain an
immediate  braid presentation of the link from the divide.

Strongly invertible links are closely related with OMS--divides. Let $L$ be an oriented link in $\S^3$ and $j$ be 
an involution of $(\S^3,L)$ with non empty fixed point set, which preserves the orientation of $\S^3$ and reverses 
the orientation of $L$. Then $(L,j)$ is called a strongly invertible link. As we said above, the link 
of a divide is strongly invertible for the implicit strong inversion $j(p,v)=(p,-v)$. Two strongly invertible 
links $(L,j)$ and $(L',j')$ are \emph{strongly equivalent} if there exists an isotopy $\varphi_t$, $t\in [0,1]$ of 
$\S^3$ sending $L$ to $L'$ such that $\varphi_1\circ j=j'\circ \varphi_1$.  

Isotopies through MS--divides give rise to strong equivalence of associated links. Also,
one can find in Couture \cite {C} an (exhaustive) list of elementary moves on MS--divides that preserve strong 
equivalence classes of the associated links. As a particular case, given a MS--divide, we can always 
construct another one, using these moves, which is an OMS--divide. Besides, we can transpose the moves on 
MS--divides directly to moves on OMS--divides. Two OMS--divides obtained one from the other by those moves on 
OMS--divides are called $\cM$--equivalent (see the list of moves in \fullref{S1:SS3}). Then as an essential result of \cite {C}, we have:
\begin{theorem*}[cf.  \cite {C}]
\begin{enumerate}
	\item Every strongly invertible link is strongly equivalent to the link 
	of an OMS--divide.
	\item The links of two OMS--divides are strongly equivalent if and only if the OMS--divides are 
	$\cM$--equivalent.  
\end{enumerate}
\end{theorem*}
As the Jones polynomial is invariant under Reidemeister moves on links diagrams, there exits a Laurent polynomial 
for an OMS--divide with integral coefficients (see Couture \cite{C}), which is invariant under $\cM$--equivalence 
and so invariant under strong equivalence of strongly invertible links. Modulo 2, this polynomial coincide with 
Jones polynomial of the link. The purpose of this paper is to interpret the polynomial of an OMS--divide as the 
graded Euler characteristics of a graded complex of $\Z_2$--vector spaces (theorem \ref{S3:P5}). Besides, 
if we call Khovanov homology of an OMS--divide the homology of this complex, then we have a stronger result:
%
\begin{theorem*} {\rm \textbf{\ref{S3:T1}}}\ Khovanov homology of OMS--divides is invariant under 
$\cM$--equiva\-len\-ce.
\end{theorem*}
%
%
\begin{corollary*} {\rm \textbf{\ref{S3:C1}}}\ Khovanov homology is an invariant for strong equivalence 
of stron\-gly invertible links.
\end{corollary*}
%
Eventually, Khovanov homology of OMS--divides is a refinement of the polynomial invariant of OMS--divides. 
%
\section{Divides and OMS--divides.}\label{S1}%
%
\subsection{Divides and links of divides.}\label{S1:SS1} %
%
A \emph{divide} of the unit $2$--disc $\D^2$ of $\R^2(\simeq \C)$ is the image $\Ga$ of 
a proper generic immersion $\gamma$: 
\begin{equation}
\gamma :(J,\partial J)\to (\D^2,\partial \D^2), \qquad 
J=\Big (\bigsqcup_{j=1}^rI_j\Big )\sqcup \Big (\bigsqcup_{j=1}^sS_j\Big ).
\label{Divide}
\end{equation}
where $I_j$ and $S_j$ are respectively copies of $[0,1]$ and $\S^1=\{z\in \C:|z|=1\}$, generic meaning 
that the only singularities of $\gamma$ are ordinary double points and $\Ga=\gamma(J)$ intersects 
$\partial \D^2$ transversally. Every $\gamma(I_j)$ (resp. $\gamma(S_j)$) is called \emph{interval} (resp.
\emph{circular}) branch.

Let $S(\D^2)=\{ (p,v)\in T\D^2 : |p|^2+|v|^2=1\}$ be the unit sphere of the tangent space 
$T\D^2\simeq\D^2\times \C$. With a divide, A'Campo \cite{AC} associates a link $\cL(\Ga)$ in $S(\D^2)$:
\begin{equation}
		\cL(\Ga)=\{ (p,v)\in S(\D^2) : p\in \Ga,\ v\in T_p\Ga\} .
		\label{Link}
\end{equation}
This link has a natural orientation induced by the two possible orientations of the branches of $\Ga$ 
and is strongly invertible for the involution $j(p,v)=(p,-v)$ of $S(\D^2)$ with axis 
${\rm Fix}(j)=\partial \D^2\times \{0\}$ (see \fullref{S1:SS3} below). Each interval 
branch of $\Ga$ leads to a strongly invertible component of $\cL(\Ga)$ and each circular branch of 
$\Ga$ to two components of $\cL(\Ga)$ interchanged by $j$.
%
\subsection{OMS--Divides.}\label{S1:SS2} %
%
%
Let $\Ga$ be a divide. Suppose there exists $(a,b)\in ]0,1[\times ]0,1[$, $a^2+b^2<1$ such that:
\begin{enumerate}
	\item $\Ga\subset \{ x+iy\in \D^2:-b<y<b\}$ and the restriction $\rho_{|\Ga}$ is a Morse function;
	\item all double points of $\Ga$ are contained in $]-a,a[\times ]-b,b[$;
	\item all maxima (resp. minima) of $\rho_{|\Ga}$, called \emph{vertical tangent points}, project on $a$ 
	(resp. $-a$). 
\end{enumerate}
Then $\Ga$ is called an \emph{ordered Morse divide}. Double points and vertical 
tangent points will be called \emph{singular points} of $\Ga$. Now let $\epsilon$ be a function that associates 
$+$ or $-$ sign with each double point of $\Ga$. Then $(\Ga,\epsilon)$ is called an \emph{ordered Morse 
signed divide} (OMS--divide) (relative to the projection $\rho(x+iy)=x$) (see {\rm\cite {CP}}). 

\begin{figure}[!ht]
		\centering
		\includegraphics[width=8.5cm]{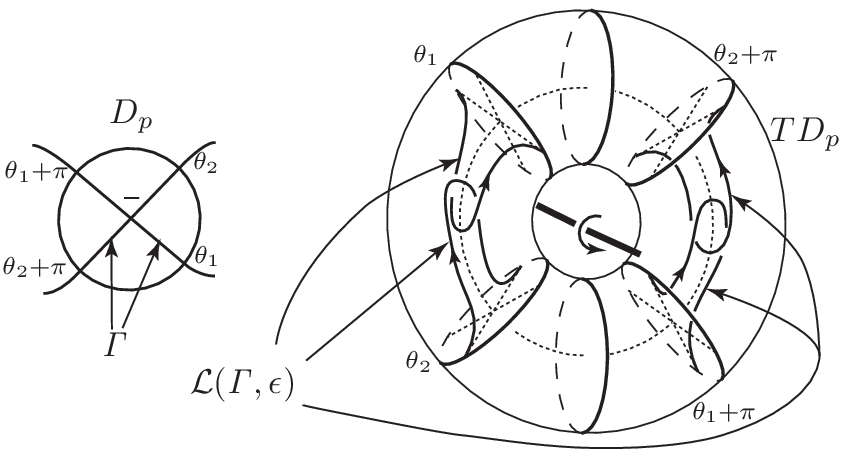}
		\caption{The link $\cL(\Ga,\epsilon)$ over a negative double point.}
		\label{fig1-2-1}
\end{figure}

Let's associate an oriented $j$--strongly invertible link $\cL(\Ga,\epsilon)$ with an OMS--divide
$(\Ga,\epsilon)$. This link coincides with $\cL(\Ga)$ except in solid tori $TD_p\cap S(\D^2)\simeq D_p\times \S^1$ 
over small discs $D_p$ around negative double points $p$ of $(\Ga,\epsilon)$ where we change the two 
$j$--symmetric crossings from over to under.
More precisely, in such a solid torus $TD_p\cap S(\D^2)$, the link is defined according to \fullref{fig1-2-1}. 
If $\epsilon=+$ for all double points then $\cL(\Ga,\epsilon)=\cL(\Ga)$.

Besides, from a divide $\Ga$ one can construct an OMS--divide 
$(\Ga',\epsilon)$ by a succession of moves and isotopies,  such that $\cL(\Ga',\epsilon)$ and $
\cL(\Ga)$ are isotopic (see \cite {CP}) by an isotopy that respects the involution $j$
(see  \fullref{S1:SS3} below for a more precise definition).
%
\rem{\label{S1:R1} For simplicity, we will only consider an OMS--divide $(\Ga,\epsilon)$ in 
$[-a,a]\times [-b,b]$, omitting its trivial part outside this rectangle. After rescaling, we also 
suppose that $a=b=1$. Since we will often consider diagrams of local parts of OMS--divides $(\Ga,\epsilon)$, 
we distinguish \emph{end points} of $\Ga$, i.e. points of $\Ga$ in $\{-1,1\}\times [-1,1]$ without vertical 
tangent by a big point (see \fullref{fig1-2-2}).
\begin{figure}[!ht]
		\centering
		\includegraphics{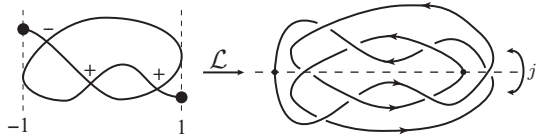}
		\caption{A representative OMS--divide for the strongly invertible knot $5_2$.}
		\label{fig1-2-2}
\end{figure}
\\
Moreover, we will simply denote $\Ga$ instead of $(\Ga,\epsilon)$ if no ambiguity 
occurs in the context.
}
%
%
\subsection{{\cal M}-equivalence for OMS--divides.}\label{S1:SS3} %
%
Two OMS--divides $\Ga$ and $\Ga'$ are \emph{$\cM$--equivalent} if we obtain one from the other by isotopy through
OMS--divides and a finite sequence of the moves described on \fullref {fig1-3-1} or symmetric situations with 
respect to horizontal and vertical directions (see \cite {C}):
\begin{figure}[!ht]
		\centering
		\includegraphics{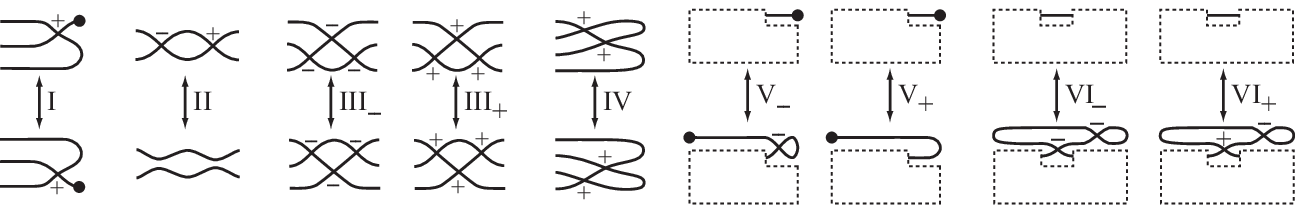}
		\caption{Moves of $\cM$--equivalence.}
		\label{fig1-3-1}
\end{figure}
Let $j$ be an orientation preserving involution of $\S^3$ with non empty fix point set \big (i.e. $\Fix (j)$ is 
trivial knot according to the solution of Smith conjecture \cite{M-B}\big ). An oriented link 
$L\subset \S^3$ is $j$--strongly invertible if $j$ sends $L$ to itself with opposite orientation. The couple 
$(L,j)$ is called a strongly invertible link. With the link of an OMS--divide, we implicitly associate natural 
orientation and involution $j(p,v)=(p,-v)$ as in \fullref {S1:SS1}: such a link is strongly invertible. 

Two strongly invertible links $(L,j)$ and $(L',j')$ are called \emph{strongly equivalent}\footnote{A 
same link $L$ may have two strong inversions $j$ and $j'$ such that $(L,j)$ and $(L,j')$ are not strongly 
equivalent.} if there exists and isotopy $\varphi_t$, $t\in [0,1]$ of $\S^3$ sending $L$ to $L'$ such that 
$\varphi_1\circ j=j'\circ \varphi_1$. One can easily prove that $\cM$--equivalent OMS--divides give rise to strongly equivalent strongly invertible links. Conversely, let's recall the following crucial theorem relating OMS--divides with strongly invertible links.
%
\thm{\label{S1:SS3:T1} {\rm \cite{C}}
\begin{enumerate}
	\item Every strongly invertible link is strongly equivalent to the link 
	of an OMS--divide.
	\item The links of two OMS--divides are strongly equivalent if and only if the OMS--divides are 
	$\cM$--equivalent.  
\end{enumerate}
}
%
\section{The polynomial of an OMS--divide.}\label{S2}%
\begin{figure}[!ht]
\centering
	\includegraphics{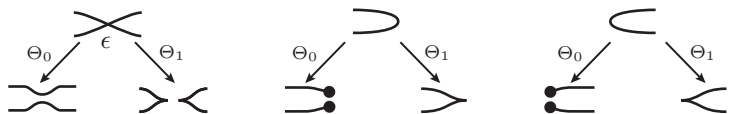}
	\caption{Local splittings.}
	\label{fig2-1}
\end{figure}

\begin{figure}[!ht]
	\centering
	\includegraphics{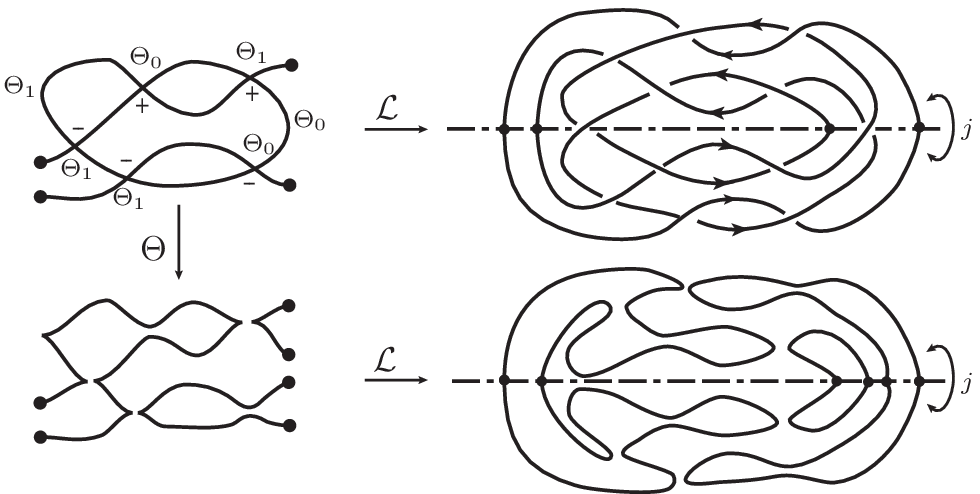}
	\caption{The links $\cL(\Ga,\epsilon)$ and $\cL(\Theta(\Ga,\epsilon))$, $\Theta=\Theta_{1110010}$.}
	\label{fig2-2}
\end{figure}

%
Let's denote by $\Theta_0$ and $\Theta_1$ the local splittings of an OMS--divide $(\Ga,\epsilon)$  in a  
neighborhood of a double point or vertical tangent point described in \fullref {fig2-1}
($\Theta_0$ ``smoothes'' the OMS--divide whereas $\Theta_1$ introduces horizontal cusps).
%
\defin{\label{S2:D1} 
\begin{enumerate}\noespace
	\item We extend the notion of OMS--divide: a \emph{cuspidal divide} $\Ga:=(\Ga,\epsilon)$ is 
	a signed diagram like an OMS--divide except that it has a finite number of horizontal 
	cusps (as in the result of type $\Theta_1$ splittings). For instance,	a partially (or totally) 
	transformed OMS--divide through $\Theta_0$ and $\Theta_1$ is a cuspidal divide. 
	\item Let $(\Ga,\epsilon)$ be an OMS--divide (or more generally of a cuspidal divide) 
	with double and vertical tangent points numbered by $p_1,\dots ,p_n$. Let $[k]$ be the word 
	$k_1k_2\dots k_n$\quad $k_i\in \{0,1\}$. A \emph{state} $(S,\Theta_{[k]})$ of $(\Ga,\epsilon)$ 
	is the combination of:\\
	$\bullet$ a succession of local splittings $\Theta_{[k]} =(\Theta_{k_1},\dots ,\Theta_{k_n})$
	at $p_1,\dots ,p_n$.\\
	$\bullet$ the cuspidal divide  $S=\Theta_{[k]}(\Ga,\epsilon)$ without double points nor vertical tangent points 
	obtained by transforming $\Ga$ through $\Theta_{[k]}$.\\ 
	For simplification, we will often identify the cuspidal divide $S$ with the state  $(S,\Theta_{[k]})$.
	We denote by $\st(\Ga,\epsilon)$ the 	set of all states of $(\Ga,\epsilon)$.
\end{enumerate}
}
%
One can define a $j$--strongly invertible link $\cL\big (\Ga,\epsilon\big )$ associated with a cuspidal divide 
$(\Ga,\epsilon)$  exactly in the same way we have done for OMS--divide. However, such a link is generally 
unoriented precisely because of the introduction of cusps. Each local splitting at a double point of 
$(\Ga,\epsilon)$ corresponds to simultaneously smoothing two symmetric crossing points of the corresponding 
representative closed braid diagram of $\cL(\Ga,\epsilon)$ (see \cite {CP}) whereas each local splitting at a 
vertical tangent point corresponds to smoothing a crossing point through the axis of the inversion $j$ (see 
\fullref{fig2-2}).
%

Let $\Ga:=(\Ga,\epsilon)$ be an OMS--divide. Let $n=n_++n_-+n_0$ be the number of singular points of $\Ga$ where
$n_+$, $n_-$ are respectively the numbers of positive and negative double points, and $n_0$ the number of vertical 
tangent points. Let's call: 
\begin{equation}
w(\Ga)=2n_+-2n_-+n_0
	\label{w}
\end{equation}
the \emph {writhe} of $\Ga$ (i.e. the writhe of the representative closed braid diagram of $\cL(\Ga,\epsilon)$ 
(see \cite {CP}) with $2n_++2n_-+n_0$ crossings obtained from $(\Ga,\epsilon)$).

For a state $S\in \st(\Ga,\epsilon)$, let $\cl(S)$ be the number of closed connected components and $\op(S)$ 
be the number of open connected components (i.e.  with two end points). Let $r_+(S)$, 
$r_-(S)$ and $r_0(S)$ be the numbers of $\Theta_1$ local splittings (\fullref {fig2-1}) for positive 
double points, negative double points and vertical tangent points respectively to obtain $S$ from 
$(\Ga,\epsilon)$. Let's set: 
\begin{equation}
i(S)=r_+(S)-r_-(S)+r_0(S)\qquad 
k(S)=w(\Ga)+2i(S)-r_0(S).\label {ik}
\end{equation}
%
\defin{\label{S2:D2} {\rm (cf. \cite {C})} The polynomial of an OMS--divide $\Ga$ (and more generally of a cuspidal 
divide) is the Laurent polynomial (of the variable $\sqrt t$) defined by:
\begin{equation}
\dis W_{\Ga}(t)=\sum_{S\in \st(\Ga,\epsilon)}(-1)^{i(S)}(\sqrt t)^{k(S)}
\left(\frac 1t +t\right )^{\cl(S)}\left(\frac 1{\sqrt t}+\sqrt t\right )^{\op(S)-1}\label {W}
\end{equation}
}
%
%
\prop  {\label {P2}{\rm (cf. \cite {C})} The polynomial of an OMS--divide is invariant under $\cM$--equiva\-len\-ce of 
OMS--divides and so is an invariant for strong equivalence of strongly invertible links.
}
%
%
\defin{\label{S2:D3} A state $S$ with $+$ or $-$ assignment to each connected component is called an 
\emph{enhanced state}, and is denoted by $\widetilde  S$. The set of enhanced states of $(\Ga,\epsilon)$ is 
denoted by $\widetilde {\st}(\Ga,\epsilon)$, and $S$ is called the \emph{underlying state} of $\widetilde  S$.
}
%
Let $\widetilde  S$ be an enhanced state with underling state $S$. 
The numbers $i(\widetilde  S):=i(S)$ and $k(\widetilde  S):=k(S)$ in (\ref{ik}) do not depend of the signs of the 
components. The subset of enhanced states $\widetilde  S$ of $(\Ga,\epsilon)$ such that $i(\widetilde  S)=i$ is 
denoted by $\widetilde {\st}_i(\Ga,\epsilon)$.

Let's denote by $\delta_\cl(\widetilde  S)$ (resp. $\delta_\op(\widetilde  S)$) the difference 
between the number of positive and negative closed (resp. open) components of $\widetilde  S$. Then we 
define the \emph{degree} $j(\widetilde  S)$ of the enhanced state $\widetilde  S$, which depends of the signs of 
the components of $S$ by:
\begin{equation}
j(\widetilde  S)=k(\widetilde  S)+2\delta_\cl(\widetilde  S)+\delta_\op(\widetilde  S).
\label {j}
\end{equation}
We can now reformulate the polynomial of an OMS--divide $\Ga$:
\begin{equation}
W_{\Ga}(t)=\frac {\sqrt t}{1+t}\sum_{\widetilde  S\in \widetilde{\st}(\Ga)}(-1)^{i(\widetilde  S)}
(\sqrt t)^{j(\widetilde  S)}
\label {W2}
\end{equation}
\rem{\label{S2:R1} $j(\widetilde  S)$ always has the same parity as half the number of end points of 
$\Ga$. We also have the inequalities:
\begin{equation}
-n_-\le i(\widetilde  S)\le n_++n_0\qquad 
2n_+-4n_-+n_0 \le k(\widetilde  S)\le 4n_+-2n_-+2n_0.
\label {encadrement i et k}
\end{equation}
}
%
\section{Categorification.}\label{S3} 										 %
%
\subsection{Complex associated with an OMS--divide.}\label{S3:SS1}
%
In this section, we define a graded complex of $\Z_2$-vector spaces\footnote{Here we choose $\Z_2$-vector spaces 
for simplification to avoid signs. We can easily generalize taking for instance $\Z$-modules or $\Q$-vector spaces} 
 associated with a divide.
We follow here Viro's approach of Khovanov homology for links \cite {V}, based on the Kauffman 
state model for the Jones polynomial: the polynomial of a divide also have been defined in \cite{C}
by state model.

\begin{figure}[!ht]
  \centering
  \includegraphics[width=\hsize]{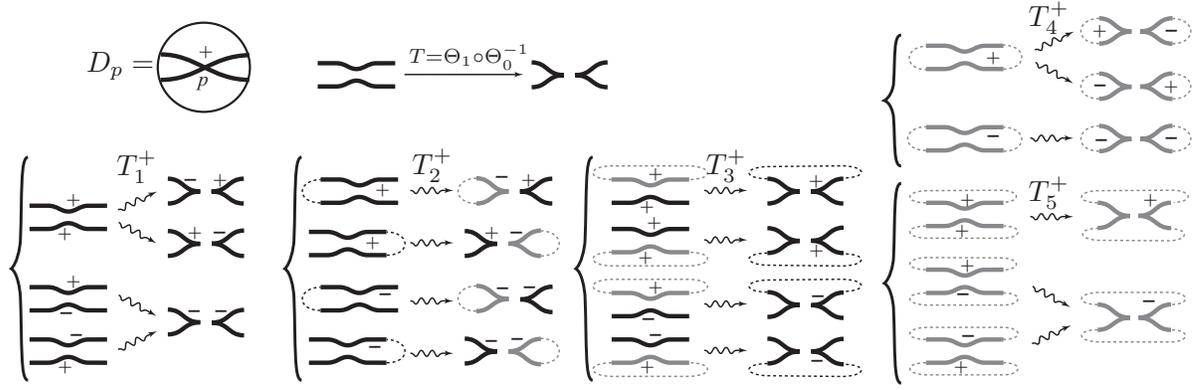}
	\caption{Case of a splitting at a positive double point $p$ ($\epsilon (p)=+$).}
	\label{fig3-1-1}
\end{figure}

\begin{figure}[!ht]
  \centering
  \includegraphics[width=\hsize]{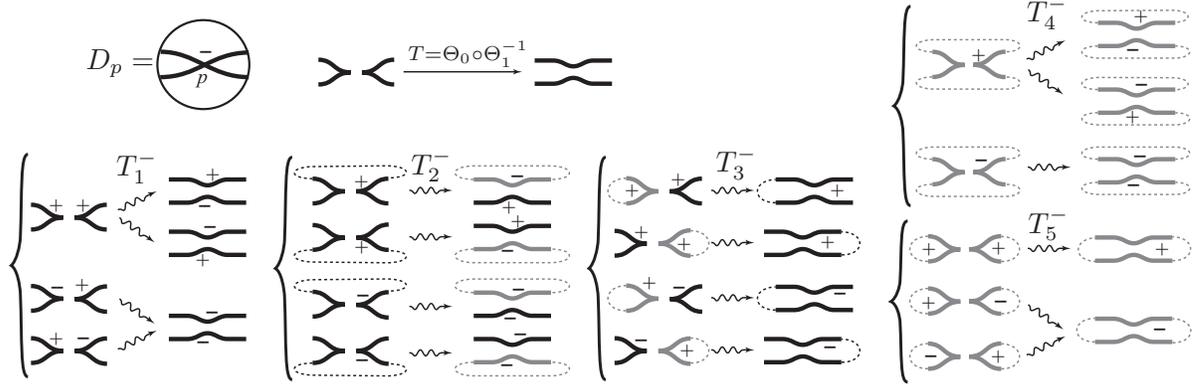}
	\caption{Case of a splitting at a negative double point $p$ ($\epsilon (p)=-$).}
	\label{fig3-1-2}
\end{figure}

\begin{figure}[!ht]
  \centering
  \includegraphics[width=\hsize]{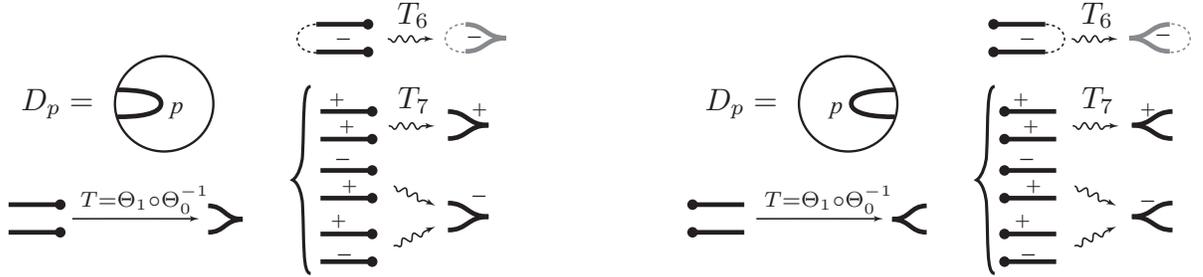}
	\caption{Case of a splitting at a vertical tangent point $p$.}
	\label{fig3-1-3}
\end{figure}

Let $\Ga:=(\Ga,\epsilon)$ be an OMS--divide (or a cuspidal divide). 
For $i\in \Z$, let $\Comp{\Ga}_i=\Z_2\{\widetilde \st_i(\Ga)\}$ be the finite dimensional 
$\Z_2$--vector space generated by enhanced states $\widetilde S$ of  $\Ga$ such that $i(\widetilde S)=i$ (if 
$i<-n_-$ or $i>n_++n_0$, $\Comp{\Ga}_i=\{0\}$). 
Degree $j(\widetilde S)$ defines a grading on $\Comp{\Ga}_i$ and we denote:
\begin{gather}
	\Comp{\Ga}=\big (\Comp{\Ga}_i\big )_{i\in \Z}\qquad \Comp{\Ga}_i=\Oplus_{j\in \Z}\Comp{\Ga}_{i,j}\\ 
	\tag*{where } \Comp{\Ga}_{i,j}=\Z_2\{\widetilde S\in \st_i(\Ga) : j(\widetilde S)=j\}.
	\label{Complex}
\end{gather}
Now we define a differential on $\Comp{\Ga}$ to obtain a (finite) complex of graded $\Z_2$--vector spaces.
%
\defin{\label {S3:D1} Let $\widetilde  S_1$,$\widetilde  S_2\in \widetilde{\st}(\Ga,\epsilon)$.  
We say that $\widetilde  S_2$ is \emph{adjacent} to $\widetilde  S_1$ ($\widetilde S_1\rightsquigarrow \widetilde 
S_2$) 
if:
\begin{enumerate}\noespace
	\item $S_1$ and $S_2$ coincide outside a neighborhood $D_p$ of a singular point $p$ of 	$(\Ga,\epsilon)$;
	\item One can pass from $\widetilde  S_1$ to $\widetilde  S_2$ by one of the following transformations $T$
		in $D_p$: 	
		\begin{itemize}
			\item[$\bullet$]  $T=\Theta_1\circ \Theta_0^{-1}$ if $p$ 
			is a positive double point or a vertical tangent point of $(\Ga,\epsilon)$;
			\item[$\bullet$]  $T=\Theta_0\circ \Theta_1^{-1}$ if $p$ 
			is a negative double point;
	\end{itemize}
	\item Signs rules described in \fullref{fig3-1-1}, \ref{fig3-1-2}, \ref{fig3-1-3} are 
fulfilled, signs of other components being unchanged.\\
{\small\noindent
(In these figures, black color is used for open components and 
gray for closed components, a dotted line means that the points are related in the state 
outside $D_p$. Lack of dotted line means that the points are not related outside $D_p$).
}
\end{enumerate}
}
%
\noindent
If $\widetilde  S_2$ is adjacent to $\widetilde  S_1$ then:
\begin{equation}
	j(\widetilde  S_1)=j(\widetilde  S_2),\qquad i(S_1)=i(S_2)-1.
	\label{RespDeg}
\end{equation}
\\
The differential $d=\big (d_i\big )_{i\in \Z}$ on $\Comp{\Ga}$, $d_i:\Comp{\Ga}_{i}\to \Comp{\Ga}_{i+1}$
is now defined in the following way: the matrix of $d_{i}$ has coefficients defined by the incidence numbers 
$(\widetilde  S_1:\widetilde  S_2)$, $\widetilde  S_1\in \widetilde \st_i(\Ga)$, 
$\widetilde  S_2\in \widetilde \st_{i+1}(\Ga)$:
\begin{equation}
	(\widetilde  S_1:\widetilde  S_2)=
								\left \{\begin{array}{ll} 
									1 &\text{if }\widetilde S_1\rightsquigarrow  \widetilde S_2\\
									0&\text{else.}
								\end{array}\right .
	\label {incidence}
\end{equation}
From (\ref{RespDeg}), $d$ respects the degree $j$, i.e. $\dis d_i=\Oplus_j d_{i,j} :\Oplus_j 
\Comp{\Ga}_{i,j}\to \Oplus_j \Comp{\Ga}_{i+1,j}$.
%
\rem{\label{S3:R1} We have dual roles for $T_1^+$ and $T_1^-$, $T^+_2$ and $T^-_3$, 
$T^+_3$ and $T^-_2$, $T^+_4$ and $T^-_5$, $T^+_5$ and $T^-_4$ in \fullref{fig3-1-1} and \fullref{fig3-1-2}. 
To go further about duality property, we could have introduced ``negative tangent points'' to interpret dual 
arrows of $T_6$ and $T_7$ in \fullref{fig3-1-3}. However, 
we didn't choose this option, since such ``negative tangent points'' can be replaced by:  
\begin{figure}[!ht]
  \centering
  \includegraphics{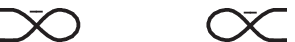}
	\label{fig3-2}
\end{figure}
\\
Also, we can see that $T_i^+$ and $T_i^-$, $i\in\{1,2,3,4\}$ give rise to analogous situations.
}
%
%
\prop {\label{S3:P1} $\dis (\Comp{\Ga},d)=\big (\Comp{\Ga}_i,d_i\big )_{i\in \Z}=
\big (\Oplus_{j\in \Z}\Comp{\Ga}_{i,j},\Oplus_{j\in \Z}d_{i,j}\big )_{i\in \Z}$ is a finite complex of 
graded $\Z_2$--vector spaces (each $\Comp{\Ga}_i$ is finitely graded by degree $j$).
}
%
\begin{proof} It suffices to verify all such diagrams:
\begin{equation}
	\begin{array}{ccc}\Ajvert {0}7 \widetilde S_1& \rightsquigarrow & \widetilde S_2 \\
										 \downsquigarrow & &\downsquigarrow \\
										 \Ajvert {15}0 \widetilde S'_2& \rightsquigarrow & \widetilde S_3
	\end{array}
	\label{CommutAdj}
\end{equation}
corresponding to splitting two singular points are commutative. Since we have $\Z_2$--vector spaces, commutative 
diagrams induce the relations $d_{i+1}d_i=0$. Notice that from the previous remark, 
we can strongly reduce the number of cases to check (see also the proof of proposition \ref{S3:P2}).  
\end{proof} 
%
\subsection{Alternative point of view.}\label{S3:SS2} 		 %
%
Here we present an alternative (more algebraic) way to see the complex $(\Comp{\Ga},d)$, in terms of 
Frobenius algebra: probably can we link complexes of OMS--divides with 1+1--TQFT (or more precisely with 
some 1+1--TQFT with symmetry property).

Let $\cA:=\Z_2\{v_-,v_+\}$ be the graded $\Z_2$--vector space generated by two elements $v_-$ and $v_+$ such that 
$\deg(v_-)=-1$ and $\deg (v_+)=1$. 
We define a commutative product $\mu_1:\cA\otimes \cA\to \cA$, a unit $\eta_1:\Z_2\to \cA$ and a non 
degenerate symmetric bilinear pairing $\beta_1:\cA\otimes \cA\to \Z_2$ by:
\begin{equation}
	\begin{array}{l}
				\mu_1(\! v_+\otimes v_+\! )=v_+, \ \
				\mu_1(\! v_+\otimes v_-\! )=\mu_1(\! v_-\otimes v_+\! )=v_-, \ \
				\mu_1(\! v_-\otimes v_-\! )=0\\
				\eta_1(0)=0,\quad
				\eta_1(1)=v_+\\
				\beta_1(\! v_+\otimes v_+\! )=0, \ \	
				\beta_1(\! v_+\otimes v_-\! )=\beta_1(\! v_-\otimes v_+\! )=1, \ \
				\beta_1(\! v_-\otimes v_-\! )=0.
	\end{array}
	\label{Prodprod1}
\end{equation}
The form $\beta_1$ induces a duality isomorphism $\cA\stackrel{\simeq}{\longleftrightarrow} \cA^*$ and $\cA$ is a 
commutative Frobenius algebra with adjoint co-product $\delta_1:\cA\to \cA\otimes \cA$ and co-unit 
$\eps_1:\cA\to \Z_2$:
\begin{equation}
		\begin{array}{l}
				\delta_1(v_+)=v_+\otimes v_-+v_-\otimes v_+	\\ 
				\delta_1(v_-)=v_-\otimes v_-
		\end{array}\quad
		\begin{array}{l}
				\eps_1(v_+)=0	\\ 
				\eps_1(v_-)=1.
		\end{array}
\label{Coprod1}
\end{equation}
Let $\phi_1 : \cA\otimes \cA\to \cA\otimes \cA$ be the flip morphism: $\phi_1(a\otimes a')= a'\otimes a$, and 
$A:\cA\to \cA$ the identity morphism.
The morphisms $\mu_1,\delta_1,\eta_1,\eps_1$  are homogeneous with respective degrees
$-1,-1,1,1$ and verify the relations of associativity, commutativity, co-associativity, co-commutativity: 
\begin{gather}
	\left\{\begin{array}{l}
		\mu_1\circ\phi_1=\mu_1\\
		\mu_1\circ (\mu_1\otimes A)=\mu_1\circ (A\otimes \mu_1)\\
		\mu_1\circ (\eta_1\otimes A)=A
	\end{array}\right .
	\
	\left\{\begin{array}{l}
		\phi_1\circ \delta_1=\delta_1\\
		(\delta_1\otimes A)\circ \delta_1=(A\otimes \delta_1)\circ \delta_1\\
		(\eps_1\otimes A)\circ \delta_1=(A\otimes \eps_1)\circ \delta_1=A
	\end{array}\right .\\
	\tag*{and}{\delta_1\circ \mu_1=(\mu_1\otimes A)\circ (A\otimes \delta_1)}.
	\label{AssComm}
\end{gather}
The vector space $\cA\otimes \cA$ has an induced structure of graded commutative 
Frobenius algebra with product, co-product, unit and co-unit~:
\begin{equation}
	\begin{array}{l}
		\mu_1^{\otimes }=(\mu_1\otimes \mu_1)\circ (A\otimes \phi\otimes A)\\
		\delta_1^{\otimes }=(A\otimes \phi\otimes A)\circ(\delta_1\otimes \delta_1)
	\end{array}\quad
	\begin{array}{l}
		\eta_1^{\otimes }=\eta_1\otimes \eta_1\\
		\eps_1^{\otimes }=\eps_1\otimes \eps_1.
	\end{array}
\label{ProdCoprod11}
\end{equation}
Let $\cB:=\Z_2\{w_-,w_+\}$ be the graded $\Z_2$--vector space generated by two elements $w_-$ and $w_+$, 
$\deg(w_-)=-2$, $\deg (w_+)=2$. Let's consider respectively the injection and surjection 
$\imath : \cB\to\cA\otimes \cA$ and $\pi:\cA\otimes \cA \to\cB$ defined by: 
\begin{equation}
\begin{array}{l}
		\imath (w_+)=v_+\otimes v_+ \\
		\imath (w_-)=v_-\otimes v_- 
\end{array}
\qquad
\begin{array}{l}
		\pi (v_+\otimes v_+)=w_+ \\
		\pi (v_+\otimes v_-)=\pi(v_-\otimes v_+)=0\\
		\pi (v_-\otimes v_-)=w_-
\end{array}
\label{InjSurj}
\end{equation}
Then $\cB$ canonically inherits from $\cA\otimes \cA$ of a structure of graded commutative 
Frobenius algebra with product, co-product, unit  and co-unit $\mu_2,\delta_2,\eta_2,\eps_2$ with
respective degrees $-2,-2,2,2$ satisfying~:
\begin{equation}
\mu_2=\pi\circ \mu_1^{\otimes }\circ (\imath\otimes \imath)\quad 
\delta_2=(\pi\otimes \pi)\circ \delta_1^{\otimes }\circ \imath\quad 
\eta_2=\pi\circ \eta_1^{\otimes }\quad 
\eps_2=\eps_1^{\otimes }\circ \imath.
\label{ProdCoprod2}
\end{equation}
The morphisms $\imath$ and $\pi$ are adjoint with degree $-2$. We denote by 
$\phi_2 : \cB\otimes \cB\to \cB\otimes \cB$ 
the flip morphism $\phi_2(b\otimes b')= b'\otimes b$, and 
$B:\cB\to \cB$ the identity morphism.

For each (non enhanced) state $S$ of $\Ga$, let's number all open components  with 
$p\in \{ 1,\dots ,\op(S)\}$ and all closed
components with $q\in \{ 1,\dots ,\cl(S)\}$. Then for an enhanced state $\widetilde S$ with 
underlying state $S$,we define the tensor product:
\begin{equation}
		t(\widetilde  S)=\Otimes_{p=1}^{\op(S)}v_{\pm_\op(p)}\otimes\Otimes _{q=1}^{\cl(S)}w_{\pm_\cl(q)}\in
		\cA^{\otimes \op(S)}\otimes\cB^{\otimes \cl(S)}
		\label {Tensor}
\end{equation}
where $\pm_\op(p)$ and $\pm_\cl(q)$ are the $+$ or $-$ signs of the $p$--th open and the $q$--th closed components
of $\widetilde  S$ respectively. The degree of $t(\widetilde  S)$ does not correspond to the degree 
$j(\widetilde  S)$
of $\widetilde  S$:
\begin{equation}
\deg (t(\widetilde  S))=\delta_\op(\widetilde  S)+2\delta_\cl(\widetilde  S)=
j(\widetilde  S)-k(\widetilde  S)=j(\widetilde  S)-k(S).
\label{Deg}
\end{equation}
So we introduce the following definition:
%
\defin {\label{S3:D2} Translation of the degree of a graded vector space.\\
Let $\cV=\Oplus\limits_{j\in \Z}\cV_j$ be a graded $\Z_2$--vector space.	The	translated 
graded $\Z_2$--vector space $\cV\{\ell\}$ is defined by: \quad 	$\cV\{\ell\}_j=\cV_{j-\ell}$.
}
%
\noindent
Now we translate the degree of $t(\widetilde  S)$ by $k(S)$ and we define:
\begin{equation}
\cC(\Ga)=\big(  \cC_i(\Ga)\big )_{i\in \Z}\quad\text{where}\quad  
\mathcal C_i(\Ga)=\Oplus_{\substack{S\in \st(\Ga) \\i(S)=i}}
\left (\cA^{\otimes \op(S)}\otimes\cB^{\otimes \cl(S)}
\right )\{k(S)\}.
\label{TensCompl}
\end{equation}
%
\prop {\label{S3:P2} The map $t: \widetilde\st(\Ga)\to \cC(\Ga)$ defined by \ref{Tensor} extends to an 
isomorphism of complexes: $t:\Comp{\Ga}\to \cC(\Ga)$.
}
%
%
\proof
The incidence relations $T^\pm_i$, $1\le i\le 5$ of \fullref{fig3-1-1}, \ref{fig3-1-2} 
induce morphisms of $\Z_2$--vector spaces denoted by $T_i$, which have degree $-2$:
\begin{equation}
	\begin{array}{l}
		T_1 = \delta_1\circ \mu_1 \quad \qquad 
		T_3 = \mu_1\circ (A\otimes \mu_1)\circ (A\otimes\imath) \quad \qquad 
		T_5 = \mu_2 \\
		\qquad T_2 = \big (A\otimes \pi)\circ (A\otimes\delta_1)\circ \delta_1  \qquad \qquad  T_4 = \delta_2.
	\end{array} 
\label{transfoT1}
\end{equation}
and the incidence relations  $T_6,T_7$  of \fullref{fig3-1-3} induce morphisms denoted by the same letters 
$T_6,T_7$, which have degree $-1$: 
\begin{equation}
	T_6 = \pi\circ \delta_1\qquad T_7 = \mu_1.  
	\label{transfoT2}
\end{equation}
More precisely, we have:
\begin{gather}
	\begin{array}{r@{}c@{\,}c@{\ }l}
		T_1 : & \cA\otimes \cA & \to     & (\cA\otimes \cA ) \\
						&\scr v_+\otimes\ v_+ & \mapsto & \scr v_+\otimes\ v_-+v_-\otimes\ v_+\\
						&\left .\begin{array}{r}
								\scr v_+\otimes\ v_- \\ 
								\scr v_-\otimes\ v_+
		  				\end{array}\hskip -6pt\right \}  
		  											 & \scr \mapsto & \scr v_-\otimes\ v_-\\
						&\scr  v_-\otimes\ v_- &\mapsto  &\scr  0
	\end{array}
	\quad 
	\begin{array}{r@{\, }c@{\,}c@{\ }l}
		T_3 :& \cA\otimes\cB  & \to     & \cA \\
						& \scr v_+\otimes\ w_+ & \mapsto & \scr v_+ \\
						& \scr v_-\otimes\ w_+ & \mapsto & \scr v_- \\
						& \left .\begin{array}{r}
								\scr v_+\otimes\ w_- \\ 
								\scr v_-\otimes\ w_-
		  				\end{array}\hskip -6pt\right \}
		  											 & \mapsto & \scr 0
	\end{array}
	\quad
	\begin{array}{l@{\, }c@{\, }c@{\ }l}
		T_5 : 	& \cB\otimes \cB & \to     & \cB \\
							& \scr w_+\otimes\ w_+ & \mapsto & \scr w_+  		\\
							&\left .\begin{array}{r}
					\scr w_+\otimes\ w_- \\ 
					\scr w_-\otimes\ w_+
			  \end{array}\hskip -6pt\right \}
			  											 & \mapsto & \scr w_-  \\
							& \scr w_-\otimes\ w_- & \mapsto & \scr 0
	\end{array}
	\nonumber\\
	\begin{array}{l@{\ }c@{\ }c@{\ }l}
		T_2 :   & \cA & \to     & (\cA\otimes \cB) \\
						& \scr v_+ & \mapsto & \scr v_+\otimes\ w_- \\
						& \scr v_- & \mapsto & \scr v_-\otimes\ w_- 				
	\end{array}
	\quad 
	\begin{array}{l@{\ }l@{\ }c@{\ }l}
		T_4 : & \cB & \to     & (\cB\otimes \cB)         \\
						& \scr w_+ & \mapsto & \scr w_+\otimes\ w_-+w_-\otimes\ w_+ \\
						& \scr w_- & \mapsto & \scr w_-\otimes\ w_-                
	\end{array}
	\\
	\begin{array}{l@{\ }c@{\ }c@{\ }l}
	T_6  :  & \cA &\to  & \cB \\
						& \scr v_+ & \mapsto & \scr 0   \\
						& \scr v_- & \mapsto & \scr w_- 
	\end{array}
	\quad
	\begin{array}{r@{\, }c@{\, }c@{\ }l}
		T_7 \, :& \cA\otimes\cA  & \to     & \cA \\
						& \scr v_+\otimes\ v_+ & \mapsto & \scr v_+ \\
						&\left .\begin{array}{r}
									\scr v_+\otimes\ v_- \\ 
									\scr v_-\otimes\ v_+
							\end{array}\hskip -6pt\right \}
														 & \mapsto & \scr v_- \\
						& \scr v_-\otimes\ v_- & \mapsto &\scr  0
	\end{array}
	\nonumber
	\label{transfoT3}
\end{gather}
Using these morphisms, we transfer the differential on $\Comp{\Ga}$ to a differential on $\cC(\Ga)$.
Notice according to remark \ref{S3:R1} that $T_1$ is self-adjoint and that $T_2$ and $T_3$ 
(resp. $T_4$ and $T_5$) are adjoint. Moreover, $T_2$ and $T_4$ are injective whereas $T_3$ and $T_5$ are 
surjective. Also $T_7$ is surjective. The relation $d\circ d=0$ \big (induced by commutative diagrams 
(\ref{CommutAdj}) in the proof of proposition \ref{S3:P1}\big ) is recovered using the following 
relations:

-- symmetry properties:
\begin{equation*}
	\begin{array}{l@{\qquad}l}
		T_1\circ \phi_1=\phi_1\circ T_1=T_1& T_4=\phi_2\circ T_4 \\
		(A\otimes T_2)\circ \phi_1=(\phi_1\otimes B)\circ (A\otimes T_2) & T_5=T_5\circ \phi_2\\
		\phi_1 \circ (A\otimes T_3)=(A\otimes T_3)\circ  (\phi_1\otimes B) & T_7=T_7\circ \phi_1
	\end{array}
\end{equation*}
-- commutativity properties corresponding to the splitting of two double points:	
\begin{gather*}
	\begin{array}{@{}>{\Ajvert {0}7}l@{\qquad}l}
	(T_1\otimes A)\circ (A\otimes T_1)=(A\otimes T_1)\circ (T_1\otimes A) &
	(T_4\otimes B)\circ T_4=(B\otimes T_4)\circ T_4\\
	(A\otimes T_2)\circ T_1=(T_1\otimes B)\circ (A\otimes T_2) &
	T_5\circ (T_5\otimes B)=T_5\circ (B\otimes T_5)\\ 
	T_1\circ (A\otimes T_3)=(A\otimes T_3)\circ (T_1\otimes B)&
	T_4\circ T_5=(B\otimes T_5)\circ (T_4\otimes B)
	\end{array}
\end{gather*}
\begin{gather*}
	T_1\circ T_1=(A\otimes T_3)\circ (\phi_1\otimes B)\circ (A\otimes T_2)=0 \\
	(T_2\otimes B)\circ T_2=(A\otimes \phi_2)\circ (T_2\otimes B)\circ T_2=(A \otimes T_4)\circ T_2\\ 
	T_3\circ	(T_3\otimes B)=T_3\circ (T_3\otimes B)\circ (A\otimes \phi_2)=T_3\circ (A \otimes T_5)\\
	\begin{array}{@{}>{\Ajvert {0}7}r@{\ }c@{\ }l}
	T_2\circ T_3&=&(T_3\otimes B)\circ (A\otimes \phi_2)\circ (T_2\otimes B)\\
			&=&(T_3\otimes B)\circ (A\otimes T_4)=(A\otimes T_5)\circ(T_2\otimes B)
	\end{array}
\end{gather*}
-- commutativity properties corresponding to the splitting of a double point and a vertical tangent point:	
\begin{equation*}
	\begin{array}{@{}>{\Ajvert {0}7}l@{\qquad}l} 
		T_4\circ T_6=(T_6\otimes B)\circ T_2 & T_7\circ (A\otimes T_3)=T_3\circ (T_7\otimes B)\\
		T_6\circ T_3=T_5\circ (T_6\otimes B) &	T_2\circ T_7=(T_7\otimes B)\circ (A\otimes T_2) \\	
		T_2\circ T_7 =(A\otimes T_6)\circ T_1 & T_7\circ T_1=0=T_3\circ (A\otimes T_6)
	\end{array}
\end{equation*}
-- commutativity properties corresponding to the splitting of two vertical tangent points:	
$$
	T_1\circ (A\otimes T_7)=(A\otimes T_7)\circ (T_1\otimes A).
\rlap{\hspace {3.4cm}\qedsymbol}
$$
%
\rem{\label{S3:R2}
The units and co-units $\eta_1,\eta_2,\eps_1,\eps_2$ of $\cA$ and $\cB$
correspond respectively to the creation of a positive open component, the creation of a positive
closed component, the destruction of a negative open component and the destruction of a 
negative closed component. Besides, $A\otimes \eps_2$ and $B\otimes \eps_2$ are left inverses of $T_2$ and $T_4$
whereas $A\otimes \eta_2$, $B\otimes \eta_2$ and $A\otimes \eta_1$ are right inverses of $T_3$, $T_5$ and $T_7$.
In \fullref{S4}, we will often refer to these morphisms together with the following ones:
$$
\begin{array}{r@{\ }c@{\ }c@{\ }c}
\wbar \eta_1 : &\Z_2 & \to   	& \cA\\
					  &\scr 1 		& \mapsto & \ \scr v_- \\
\end{array}
\quad
\begin{array}{r@{\ }c@{\ }c@{\ }c}
\wbar \eta_2 : & \Z_2 & \to     & \cB\\
								 & \scr 1    & \mapsto & \ \scr w_- \\
\end{array}
\quad 
\begin{array}{r@{\ }c@{\ }c@{\ }c}
\wbar \varepsilon_1 : & \cA& \to     & \Z_2 \\
					      			  &\scr v_+ 	    & \mapsto & \scr 1    \\
					        			&\scr v_- 	    & \mapsto & \scr 0
\end{array}
\quad
\begin{array}{r@{\ }c@{\ }c@{\ }c}
\wbar \varepsilon_2 : & \cB& \to     & \Z_2 \\
									& \scr w_+ 		 & \mapsto & \scr 1    \\
									& \scr w_- 		 & \mapsto & \scr 0    
\end{array}
$$ 
which correspond respectively to the creation of a negative open component, the creation of a negative 
closed component, the destruction of a positive open component and the destruction of a 
positive closed component, and to the following composed morphisms: 
$$
\begin{array}{r@{\ }c@{\ }c@{\ }c}
\tau=\eta_1\eps_1 : & \cA & \to   	& \cA \\
					 						& \scr v_+ & \mapsto & \scr 0   \\
					 						& \scr v_- & \mapsto & \ \scr v_+
\end{array}
\quad
\begin{array}{r@{\ }c@{\ }c@{\ }c}
\sigma=\wbar \eta_1\eps_2 : & \cB  & \to     & \cA \\
														  & \scr w_+  & \mapsto & \scr 0   \\
														  & \scr w_-  & \mapsto & \ \scr v_-
\end{array}
$$
}
%
%
\subsection{Review of basic facts about complexes.}\label {S3:SS3} %
Let $\cC:=(\cC,d)=(\cC_i,d_i)_{i\in \Z}$ be a complex of $\Z_2$--vector spaces. We denote by $\cH(\cC)$ 
its homology:
\begin{equation}
	\cH(\cC)=(H_i)_{i\in \Z} \qquad  H_i=\Ker d_i/\Im d_{i-1}.
	\label{Homol}
\end{equation}
A complex is \emph{acyclic} if its homology is null.
%
\defin{\label {S3:D3} Shift of the grading of a complex.\\
Let $\dis (\cC,d)=(\mathcal C_i,d_i)_{i\in \Z}$ be a complex of $\Z_2$--vector spaces.  We define
the complex:
$$
(\cC,d)[k]=(\cC[k],d[k])\quad \text{by} \quad \cC[k]_i=\cC_{i-k}\quad \text{and}\quad d[k]_i=d_{i-k}.
$$	
({\small If $\dis (\cC,d)=\big (\Oplus_{j\in \Z}\mathcal C_{i,j},\Oplus\limits _{j\in \Z}d_{i,j}\big )_{i\in \Z}$
is a complex of graded $\Z_2$--vector spaces, then we can translate twice the grading of the complex and 
the degree of the vector spaces :
$$
	(\cC,d)[k]\{\ell\}=(\cC,d)\{\ell\}[k]\quad \text{is defined by}\quad
	\cC[k]\{\ell\}_{i,j}=\cC_{i-k,j-\ell}\quad
	d[k]\{\ell\}_{i,j}=d_{i-k,j-\ell}\text{ )}.
$$
}
}
%
A morphism of complexes of $\Z_2$--vector spaces $f:(\cC^0,d^0)\to (\cC^1,d^1)$ is a sequence $f=(f_i)_{i\in \Z}$ 
of linear maps \ $f_i:\cC^0_i\to \cC^1_i$\ such that\footnote{Since we are working with $\Z_2$
field, commutativity and anti-commutativity coincide so that we have equivalently \ $d^1f+fd^0=0$.}:\ 
$f\, d^0=d^1\, f$ (i.e. $\forall i$, $d^1_i\, f_i=f_{i+1}\, d^0_i$). 
%
\defin{\label {S3:D4} For a morphism of complexes $f:(\cC^0,d^0)\to (\cC^1,d^1)$, the \emph {cone} of $f$ 
is the complex denoted by $\cone(f)=(\cC_i,D_i)_{i\in \Z}$ and defined by~:
\begin{equation}
	\cC_i=\cC^0_i\oplus \cC^1_{i-1}=\cC^0_i\oplus (\cC^1[1])_{i}\qquad 
	D_i=\left(\begin{array}{cc}d_i^0&0\\f_i&d_{i-1}^1\end{array}\right ).
	\label{cone}
\end{equation}
(Notice that $(\cC^0,d^0)$ and $(\cC^1,d^1)[1]$ are sub-complexes of $(\cC,D)$).
}
%
A morphism of complexes $f:(\cC^0,d^0)\to (\cC^1,d^1)$ induces an isomorphism in homology if and only if 
$\cone(f)$ is acyclic. This is the case if $f$ is a \emph{homotopy equivalence} of complexes, i.e. there exist 
a morphism of complexes $g:(\cC^1,d^1)\to (\cC^0,d^0)$ and sequences $h^0=(h_i^0)_{i\in \Z},h^1=(h_i^1)_{i\in \Z}$ 
of linear maps (homotopies) 
$h_i^0:\cC_{i+1}^0\to \cC_i^0$ and 
$h_i^1:\cC_{i+1}^1\to \cC_i^1$ such that:
\begin{gather}
	gf=\id+h^0d^0+d^0h^0\qquad \text{and}\qquad fg=\id+h^1d^1+d^1h^1.\\
	({\scriptstyle \text{i.e.}\quad \forall i\quad g_if_i=\id+h_i^0d_i^0+d_{i-1}^0h_{i-1}^0\quad \text{and}\quad 
		f_ig_i=\id+h_i^1d_i^1+d_{i-1}^1h_{i-1}^1})\nonumber
	\label{HomotEquiv} 
\end{gather}
%
\rem{\label{S3:R3} As a particular case, if $h^0=0$, the complex $(\cC^0,d^0)$ is called a \emph{strong deformation 
retract} of 
$(\cC^1,d^1)$, with inclusion map $f$, retraction $g$ and homotopy map $h^1$. Besides, up to changing $h^1$ to a 
new homotopy $h$, we can always suppose that $hh=0$, $hf=0$ and $gh=0$. We will assume these properties
are always satisfied in the definition of strong deformation retraction.
}
%
%
\prop{\label {S3:P3} Let $(\wbar \cC^1,\wbar d^1)$ be a strong deformation retract of $(\cC^1,d^1)$ with 
retraction $r$, inclusion $j$ and homotopy map $h$ such that $hh=0$, $rh=0$, $hj=0$. Let 
$f:(\cC^0,d^0)\to (\cC^1,d^1)$ be a morphism of complexes.
Then $\cone(rf)$ is a strong deformation retract of $\cone(f)$ with retraction 
$R=\left (\begin{array}{cc}\id&0\\ 0&r\end{array}\right )$, inclusion 
$J=\left (\begin{array}{cc}\id&0\\ hf&j\end{array}\right )$ and homotopy 
$H=\left (\begin{array}{cc}0&0\\ 0&h\end{array}\right )$ such that $HH=0$,
$RH=0$ and $HJ=0$.
}
%
\begin{proof} Immediate.
\end{proof}
%
A \emph{double complex} $(\cC,d,\partial)$ is a sequence of complexes $(\cC^k,d^k)_{k\in \Z}$ of $\Z_2$--vector 
spaces and morphisms of complexes $(\partial ^k)_{k\in \Z}$:
	$$
	\cdots\mathop{\longrightarrow}^{\partial^{k-1}}  (\cC^{k},d^{k})\mathop{\longrightarrow}^{\partial^{k}} 
	(\cC^{k+1},d^{k+1})\mathop{\longrightarrow}^{\partial ^{k+1}}
	(\cC^{k+2},d^{k+2})\mathop{\longrightarrow}^{\partial^{k+2}}\cdots 
	$$
such that for all $k\in \Z$,\ $\partial^{k+1}\partial^k=0$. A \emph{morphism of two double complexes} 
is a sequence of morphisms of complexes $f=(f^k)_{k\in \Z}$:
	$$
	\begin{array}{c@{\ }c@{\ }c@{\ }c@{\ }c@{\ }c@{\ }c@{\ }c@{\ }c}
	\cdots
	&\displaystyle \mathop{\longrightarrow}^{\wbar\partial ^{k-1}} 
	&(\wbar\cC^k,\wbar d^k)
	&\displaystyle \mathop{\longrightarrow}^{\wbar\partial ^k} 
	&(\wbar\cC^{k+1},\wbar d^{k+1})
	&\displaystyle \mathop{\longrightarrow}^{\wbar\partial ^{k+1}}
	&(\wbar\cC^{k+2},\wbar d^{k+2})
	&\displaystyle\mathop{\longrightarrow}^{\wbar\partial ^{k+2}}
	&\cdots \\
	&&\ \bigg \downarrow \scriptstyle f^k&&\ \bigg \downarrow \scriptstyle f^{k+1}&&
	\ \bigg \downarrow \scriptstyle f^{k+2}\\
	\cdots
	&\displaystyle \mathop{\longrightarrow}^{\partial ^{k-1}} 
	&(\cC^k,d^k)
	&\displaystyle\mathop{\longrightarrow}^{\partial^k} 
	&(\wbar \cC^{k+1},d^{k+1})
	&\displaystyle\mathop{\longrightarrow}^{\partial^{k+1}}
	&(\cC^{k+2},d^{k+2})
	&\displaystyle\mathop{\longrightarrow}^{\partial^{n-1}}
	&\cdots
	\end{array}
	$$
such that for all $k\in \Z$,\ $f^{k+1}\wbar \partial ^k=\partial^kf^k$.
We also have notions of homotopy equivalence and strong deformation retraction for double complex. 
A morphism of  double complexes  $f:(\wbar\cC,\wbar d,\wbar\partial)\to (\cC,d,\partial)$ is a 
\emph{homotopy equivalence (of double complexes)} if there exists a morphism of double complexes 
$g:(\cC,d,\partial)\to (\wbar \cC,\wbar d,\wbar \partial)$ and homotopy 
maps $\wbar h=(\wbar h^k)_{k\in \Z}$, $h=(h^k)_{k\in \Z}$(sequences of morphisms of complexes) 
$\wbar h^k:(\wbar \cC^{k+1},\wbar d^{k+1})\to (\wbar \cC^k,\wbar d^k)$ and 
$\wbar h^k:(\cC^{k+1},d^{k+1})\to (\cC^k,d^k)$  such that for all $k$:
$$
g^kf^k=\id+\wbar h^k\wbar \partial ^k+\wbar \partial ^{k-1}\wbar h^{k-1}\qquad \text{and}\qquad 
f^kg^k=\id+h^k\partial ^k+\partial ^{k-1}h^{k-1}.
$$ 
If $\wbar h=0$, $(\wbar \cC,\wbar d,\wbar \partial)$ is called a \emph{strong deformation retract}
of the double complex $(\cC,d,\partial)$ with inclusion $f$ and retraction $g$.
Again, up to changing the homotopy $h$, we assume that it satisfies: $hh=0$, 
$hf=0$ and $gh=0$. 

Now we extend the definition of cone to a finite sequence of morphisms of complexes.
A double complex $(\cC,d,\partial)$ is $\partial $--finite if $(\cC^k,d^k)$ is trivial except
for a finite number of values of $k$.
%
\defin{\label {S3:D5} Let $\displaystyle (\cC^{0},d^{0})\mathop{\longrightarrow}^{\partial^{0}} \cdots 
\mathop{\longrightarrow}^{\partial ^{n-1}}(\cC^{n},d^{n})$ be a $\partial$--finite double complex.
Let's denote:
	$$
	\begin{array}{crcl}
	\wtilde\partial ^0_i:&\cC^0_i&\to &\cC^1_i\oplus \cC^2_{i-1}\oplus \cdots \oplus \cC^n_{i-n+1}\\
	&u&\mapsto&(\partial _i^0(u),0,\dots,0)
	\end{array}.
	$$
Then the \emph{cone} of $(\partial ^0,\dots ,\partial ^{n-1})$ is the complex defined by:
\begin{equation}
		\cone(\partial ^0,\dots ,\partial ^{n-1})=\cone(\wtilde\partial ^0,\cone(\partial ^1,\dots ,
	\partial ^{n-1})).
	\label{cone-gene}
\end{equation}
}
%
Suppose that $f=(f^k)_{0\le k\le n}$ is a morphism from a $\wbar \partial $--finite double complexes 
$(\wbar \cC,\wbar d,\wbar \partial)$ to a $\partial $--finite  $(\cC,d,\partial)$:
$$
(\wbar \cC^{0},\wbar d^{0})\mathop{\longrightarrow}^{\wbar \partial^{0}} \cdots 
	\mathop{\longrightarrow}^{\wbar \partial ^{n-1}}(\wbar \cC^{n},\wbar d^{n})\quad \text{ and }\quad 
	(\cC^{0},d^{0})\mathop{\longrightarrow}^{\partial^{0}} \cdots 
	\mathop{\longrightarrow}^{\partial ^{n-1}}(\cC^{n},d^{n}).
$$
Let's set $F_i=f_i^0\oplus f_{i-1}^1\oplus \cdots \oplus f_{i-n}^n$. Then $f$ induces a morphism of complexes:
$$
C(f)=(F_i)_{i\in \Z}:\cone(\wbar \partial ^0, \dots,\wbar \partial ^{n-1})\longrightarrow 
\cone(\partial ^0,\dots,\partial ^{n-1}).
$$ 
If $f^k$ are isomorphisms, $C(f)$ is also an isomorphism.	
%
\prop{\label {S3:P4} If $(\wbar \cC,\wbar d,\wbar \partial)$ is a $\wbar\partial$--finite double 
complex, $(\cC,d,\partial)$ a $\partial$--finite double complex and 
$f:(\wbar \cC,\wbar d,\wbar \partial)\to (\cC,d,\partial)$ a homotopy equivalence
with inverse $g$, then:
$$
C(f):\cone(\wbar \partial ^0, \dots,\wbar \partial ^{n-1})\longrightarrow 
\cone(\partial ^0, \dots,\partial ^{n-1})
$$
is a homotopy equivalence of complexes with inverse $C(g)$. So $C(f)$ induces an isomorphism in homology.
\\
Furthermore if $(\wbar \cC,\wbar d,\wbar \partial)$ is a strong deformation retract of $(\cC,d,\partial)$ 
with inclusion map $f$ and retraction $g$ then $\cone(\wbar \partial ^0, \dots,\wbar \partial ^{n-1})$ is a 
strong deformation retract of $\cone(\partial ^0, \dots,\partial ^{n-1})$ with inclusion map $C(f)$ 
and retraction $C(g)$, and so $C(f)$ induces an isomorphism in homology.
}
%
\proof Let $\wbar h=(\wbar h^k)_{1\le k< n}$ and $h=(h^k)_{0\le k< n}$ be homotopies associated 
with $f$ and $g$. Then we have~:
$$
\begin{array}{l!{\hspace{10pt}}l}
\wbar h_{i+1}^{k-1}\wbar d_i^{k}=\wbar d_i^{k-1}\wbar h_i^{k-1} & \wbar h_i^{k}\wbar \partial _i^{k}+
\wbar \partial _i^{k-1}\wbar h_i^{k-1}+\id=g_i^kf_i^k\\
h_{i+1}^{k-1}d_i^{k}=d_i^{k-1}h_i^{k-1} & h_i^{k}\partial _i^{k}+\partial _i^{k-1}h_i^{k-1}+\id=f_i^kg_i^k
\end{array}
$$
Let $H=(H_i)_{i\in \Z}$ be the sequence of linear maps defined by:
$$
\begin{array}{rcl}
	H_i:\cC^0_{i+1}\oplus \cC^1_{i}\oplus\cdots \oplus \cC^n_{i-n+1}&\longrightarrow
	& \cC^0_{i}\oplus \cC^1_{i-1}\oplus\cdots \oplus \cC^n_{i-n}\\
	(x_0,x_1,\dots,x_n)\qquad &\longmapsto &(h^0_i(x_1),h^1_{i-1}(x_2),\dots,h^{n-1}_{i-n+1}(x_n),0)
\end{array}
$$
and $\wbar H=(\wbar H_i)_{i\in \Z}$ defined analogously on $(\wbar \cC,\wbar d,\wbar \partial )$.  Then
if $\wbar D$ and $D$ are the differentials of $\cone(\wbar \partial ^0, \dots,\wbar \partial ^{n-1})$ and 
$\cone(\partial ^0, \dots,\partial ^{n-1})$, we have:
$$
FG=\id+HD+DH\quad \text{and}\quad GF=\id+\wbar H\wbar D+\wbar D\wbar H.\rlap{\hspace{2cm}\qedsymbol}
$$
\subsection{Fundamental splitting lemmas.}\label {S3:SS4} %
Let $(\Ga,\epsilon)$ be an OMS--divide or a cuspidal divide. Let $p$ be a double point or a vertical 
tangent point of $\Ga$. Let $\Ga^0$ and $\Ga^1$ be the cuspidal divides obtained from $\Ga$ by 
applying $\Theta_0$ and $\Theta_1$ at $p$ respectively. Then each enhanced state of $\Ga$ can be identified 
with either an enhanced state of $\Ga^0$ or $\Ga^1$ i.e.:
\begin{equation}
	\widetilde\st(\Ga)\mathop{\simeq}^{1-1}\ \widetilde \st(\Ga^0)\sqcup \widetilde{\st}(\Ga^1).
	\label{State-decomp}
\end{equation}
Consequently, $\Comp{\Ga^0}$ and $\Comp{\Ga^1}$ can be seen as sub-complexes of 
$\Comp{\Ga}$ up to translations of the grading $i$ and the degree $j$. 
More precisely, we have :
%
\lem{\label {S3:L1} Let $d$ be the differential of $\Comp{\Ga}$.
\begin{enumerate}
	\item If $p$ is a positive double point then  $d$ 
	induces the differentials $d^0$ and $d^1$ of $\Comp{\Ga^0}\{2\}$ and $\Comp{\Ga^1}\{4\}$
	and a morphism $\dis \Comp {\Ga^0}\{2\}\mathop{\longrightarrow}^{d^{\bullet}}\Comp {\Ga^1}\{4\}$	
	such that $\Comp{\Ga}=\cone(d^{\bullet})$.
	\item If $p$ is a negative double point then  $d$ 
	induces the differentials $d^0$ and $d^1$ of $\Comp{\Ga^0}\{-2\}$ and $\Comp{\Ga^1}\{-4\}$
	and a morphism $\dis \Comp {\Ga^1}\{-4\}\mathop{\longrightarrow}^{d^{\bullet}}\Comp {\Ga^0}\{-2\}$
	such that $\Comp{\Ga}=\cone(d^{\bullet})[-1]$.
	\item If $p$ is vertical tangent point, then  $d$ 
	induces the differentials $d^0$ and $d^1$ of $\Comp{\Ga^0}\{1\}$ and $\Comp{\Ga^1}\{2\}$
	and a morphism $\dis \Comp {\Ga^0}\{1\}\mathop{\longrightarrow}^{d^{\bullet}}\Comp {\Ga^1}\{2\}$
	such that $\Comp{\Ga}=\cone(d^{\bullet})$.
\end{enumerate}
}
%
\begin{proof}
Suppose that $p$ is a positive double point of $\Ga$. Then $\Ga^0$ and $\Ga^1$ have one positive 
double point less than $\Ga$ so that the writhes of $\Ga$, $\Ga^0$ and $\Ga^1$ are related by:
$$
w(\Ga)=w(\Ga^0)+2=w(\Ga^1)+2
$$
Let $\widetilde S$ be an enhanced state of $\widetilde\st_i(\Ga)$ with degree $j=j(\widetilde S)$. 
If $\widetilde S$ is obtained from $\Ga$ using $\Theta_0$ (resp. $\Theta_1$) at $p$, then $\widetilde S$ 
can be seen as an enhanced state of $\widetilde\st_i(\Ga^0)$ with degree $j-2$
(resp. of $\widetilde\st_{i-1}(\Ga^1)$ with degree $j-4$). Besides, if $\widetilde S_0\in \widetilde\st_i(\Ga)$ 
and $\widetilde S_1\in \widetilde \st_{i+1}(\Ga)$ are adjacent enhanced states of degrees $j$ then it 
involves three cases:
\begin{itemize}
	\item either $\widetilde {S}_0\in \widetilde \st_i(\Ga^0)$ and $\widetilde {S}_1\in	\widetilde{\st}_{i+1}(\Ga^0)$ 
	are adjacent enhanced states of $\Ga^0$ with degrees $j-2$, so the differential $d^0$ of
	$\Comp{\Ga^0}\{2\}$ coincide with the restriction of $d$ to $\Comp{\Ga^0}\{2\}$;
	\item either $\widetilde S_0\in \widetilde\st_{i-1}(\Ga^1)$ and $\widetilde S_1\in \widetilde\st_i(\Ga^1)$ 
	are adjacent enhanced states of $\Ga^1$ with degrees $j-4$, so the differential $d^1$ of
	$\Comp{\Ga^1}\{4\}[1]$ coincide with the restriction of $d$ to $\Comp{\Ga}\{4\}[1]$;
	\item or $\widetilde S_0\in \widetilde\st_i(\Ga^0)$ with degree $j-2$ and $\widetilde S_1\in 	
	\widetilde\st_i(\Ga^1)$ with degree $j-4$, then $d$ induces a map $d^{\bullet}:\Comp{\Ga^0}\{2\}\to 
	\Comp{\Ga^1}\{4\}$ which is a morphism of complexes since from the proof of proposition \ref{S3:P1} $dd=0$ 	
	implies $d^{\bullet}d^0=d^1d^{\bullet}$.
\end{itemize}
Hence $\Comp{\Ga}=\cone(d^{\bullet})$.
Similar arguments hold for the two other cases.
\end{proof}
More generally, consider $k=k_++k_-+k_0$ double vertical tangent points $p_1,\dots ,p_k$
such that the $k_+$ first ones are positive double points, the next $k_-$ ones negative double points 
and the last $k_0$ ones vertical tangent points. For each words $[a]=a_1a_2\dots a_{k_+}$, $a_i\in \{0,1\}$, 
$[b]=b_1b_2\dots b_{k_-}$, $b_i\in \{0,1\}$, $[c]=c_1c_2\dots c_{k_0}$, $c_i\in \{0,1\}$, let 
$[a][b][c]$ be the word obtained by concatenation of $[a],[b],[c]$ and denote by 
$(\Ga^{[a][b][c]},\epsilon^{[a][b][c]})$ the cuspidal divide obtained from $(\Ga,\epsilon)$ by performing:
$$
\left\{ 
\begin{array}{lll}
	\Theta_{a_i} &\text{splitting at } p_i & \text{for }1\le i\le k_+ \\
	\Theta_{b_i} &\text{splitting at } p_i & \text{for }k_+< i\le k_++k_- \\
	\Theta_{c_i} &\text{splitting at } p_i & \text{for }k_++k_-< i\le k=k_++k_-+k_0.
\end{array}
\right.
$$
Let $\dis 1_{[a]}$, $1_{[b]}$ and $1_{[c]}$ be the numbers of occurrences of 1 in $[a]$, $[b]$ and $[c]$ and 
$\text{gr}([a][b][c])=1_{[a]}-1_{[b]}+1_{[c]}$.
By restriction, the differential $d$ of $\Comp{\Ga}$ coincide with the differential $d^{[a][b][c]}$ of 
$\Comp{\Ga^{[a][b][c]}}$. By iterating lemma \ref{S3:L1}, using same arguments, just following the incidence 
relations, we have:
%
\lem{\label {S3:L2} For each $\ell$, $-k_-\le \ell\le k_++k_0$, we can identify the complex:
	$$(\cC^\ell,D^\ell)=
	\Oplus_{\text{gr}([a][b][c])=\ell}\Comp{\Ga^{[a][b][c]}}\big \{2(1_{[a]}-1_{[b]}+k_+-k_-)+1_{[c]}+k_0\big \}
	$$
with a sub-complex of $\Comp{\Ga}$, with differential 	
$\dis D^{\ell}=\Oplus_{\text{gr}([a][b][c])=\ell}d^{[a][b][c]}$. 
The differential $d$ induces a structure of double complex : 
	$$
	\big (\cC^{-k_-},D^{-k_-}\big )\mathop{\rightarrowvar {.7cm}}^{\Delta^{-k_-}}
	\big (\cC^{-k_-+1},D^{-k_-+1}\big )\mathop{\rightarrowvar {1.1cm}}^{\Delta^{-k_-+1}}\ \cdots\ 
	\mathop{\rightarrowvar {1.1cm}}^{\Delta^{k_++k_0-1}}\big (\cC^{k_++k_0},D^{k_++k_0}\big )
	$$
such that $\Comp{\Ga}=\cone\big (\Delta^{-k_-},\dots,\Delta^{k_++k_0-1}\big )[-k_-]$. 
}
%
In the sequel, such a double complex will be called a \emph{splitting diagram} of 
$\Comp{\Ga}$.
%
\subsection{Khovanov homology of OMS--divides.}\label{S3:SS5} %
%
\defin {\label{S3:D6} We call \emph{Khovanov homology} $\cH(\Ga)$ of an OMS--divide (or a cuspidal divide)
$\Ga=(\Ga,\epsilon)$ the 
homology of the complex $\Comp{\Ga}=\big (\Comp{\Ga}_i\big )_{i\in \Z}$:
\begin{equation}
	\cH(\Ga)=\big (\cH_i(\Ga)\big)_{i\in \Z}\quad 
	\cH_i(\Ga)=\Oplus_{j\in \Z}\cH_{i,j}(\Ga)\quad \cH_{i,j}(\Ga)=
	\Ker d_{i,j}/\Im d_{i-1,j}.
	\label {HomoKhov}
\end{equation}
}
%
%
\prop {\label{S3:P5} If $\Ga=(\Ga,\epsilon)$ is an OMS--divide, then the  polynomial $W_{\Ga}$ and 
the graded Euler characteristics of $\cH(\Ga)$ are related by:
\begin{equation}
W_{\Ga}(t^2)=\frac {t}{1+t^2}\ \chi_{gr}(\cH(\Ga))=\frac {t}{1+t^2}\ 
\sum_{i\in \Z}(-1)^i\dim _{gr}\cH_i(\Ga)
\label{KhovPol}
\end{equation}
where the graded dimension is:\quad 
$\dis
\dim _{gr}\cH_i(\Ga)=\sum _{j\in \Z}t^j\dim _{\Z_2}\cH_{i,j}(\Ga).
$
}
%
\begin{proof} Immediate from formula (\ref {W2}).
\end{proof}
We can now formulate our main theorem:
%
\thm {\label{S3:T1} Khovanov homology of OMS--divides is invariant under $\cM$--equi\-va\-len\-ce.
}
%
Combined with theorem \ref{S1:SS3:T1}, we obtain:
%
\cor {\label{S3:C1} Khovanov homology of OMS--divides is an invariant for strong equivalence of 
strongly invertible links.
}
%
\fullref{S4} is devoted to the proof of this theorem. Notice that from proposition \ref{S3:P5}, 
this theorem \ref{S3:T1} is a refinement of proposition \ref{P2}.
%
\subsection{Examples.}\label{S3:SS6} %
\begin{figure}[!ht]
	\centering
  \includegraphics{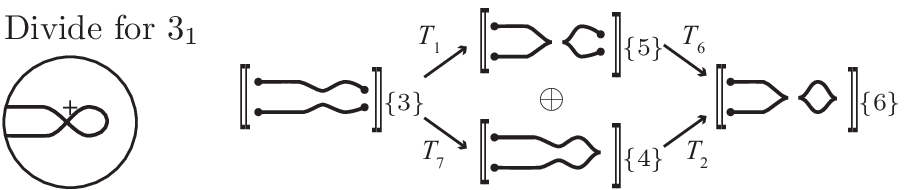}
	\caption{}
	\label{fig3-3} 
\end{figure}

\begin{figure}[!ht]
  \centering
  \includegraphics{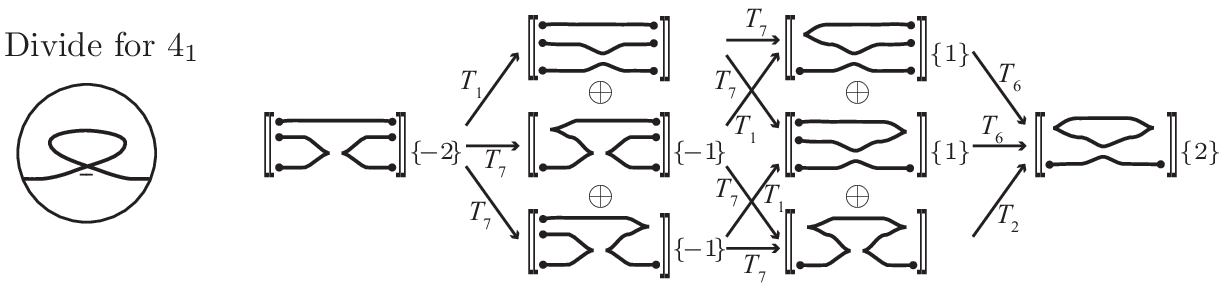}
	\caption{}
	\label{fig3-4}
\end{figure}

(1) \fullref{fig3-3} shows a divide for the link $3_1$, and its splitting diagram.
$$
\begin{array}{l}
\text{The associated complex and homology entries are:}\\
\\
\dis (\cA\otimes\cA)\{3\}\mathop{\to}^{d^0}(\cA\otimes\cA)\{5\}\oplus \cA\{4\}
\mathop{\to}^{d^1}(\cA\otimes\cB)\{6\}\\
\\
\end{array}
\ 
\begin{array}{|c|c|c|c|c|c|}
\hline
\scr \quad i \ \backslash\  j\! &\scr 1 &\scr 3 &\scr 5 &\scr 7 &\scr 9 \\
\hline
\scr 0 &\scr \Z_2 &\scr \Z_2 &  &  &  \\
\hline
\scr 1 &  &  & \scr \Z_2 &\scr \Z_2 &  \\
\hline
\scr 2 &  &  &  &\scr \Z_2 &\scr \Z_2 \\
\hline
\end{array}  
$$
(2) \fullref{fig3-4} shows a divide for the link $4_1$ and its splitting diagram.
The associated complex is:
{\small
$$
\dis (\cA^{\otimes3})\{-2\}\mathop{\to}^{d^{-1}}
(\cA^{\otimes3})\oplus(\cA^{\otimes2}\oplus\cA^{\otimes2})\{-1\}\mathop{\to}^{d^0}
(\cA^{\otimes2}\oplus\cA^{\otimes2})\{1\}\oplus \cA\mathop{\to}^{d^1}
(\cA\otimes\cB)\{2\}
$$
$$ 
\llap{and homology entries: }\quad
\begin{array}{|c|c|c|c|c|c|c|}
\hline
\quad \scr i \ \backslash\  j\! &\scr  -5 &\scr  -3 &\scr  -1 &\scr  1 &\scr  3 &\scr  5 \\
\hline
\scr  -1 &\scr  \Z_2\ &\scr \Z_2\ &  &  &  &  \\
\hline
\scr  0 &  &\scr  \Z_2 & \!\scr  (\Z_2)^2\! \! &\scr  \Z_2 &  &  \\
\hline
\scr  1 &  &  &  & \scr  \Z_2\ &\scr  \Z_2 &  \\
\hline
\scr  2 &  &  &  &  &\scr  \Z_2\ & \scr  \Z_2\  \\
\hline
\end{array}  
$$}
\section{Invariance under $\cM$--equivalence.}\label{S4}    %
\subsection{Invariance under type I moves.}\label {S4:SSI} %
\begin{figure}[!ht]
  \centering
  \includegraphics{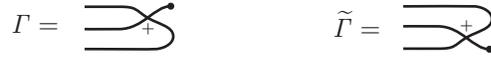}
	\caption{Type I move.}
	\label{fig-I-1}
\end{figure}

\begin{figure}[ht]
\centering
  \includegraphics{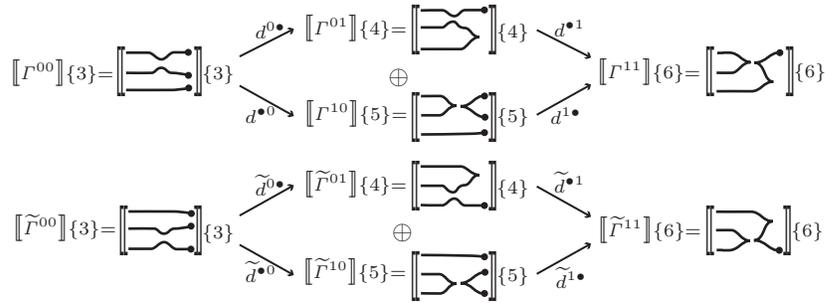}
	\caption{Splitting diagram for type I move.}
	\label{fig-I-2}
\end{figure}

\begin{figure}[!ht]
  \centering
  \includegraphics[width=\hsize]{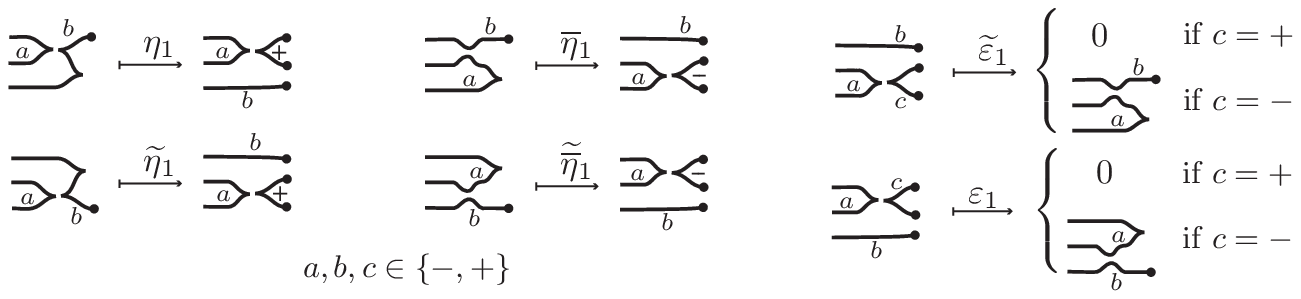}
	\caption{}
	\label{fig-I-3}
\end{figure}

Let $\Ga$ and $\wtilde\Ga$ be OMS--divides which differ only by a type I move (see \fullref{fig-I-1}).
%
\prop{\label{S4:SSI:P1}  The complexes $\Comp{\Ga}$ and $\Comp{\wtilde{\Ga}}$ have 
the same homology.
}
%
Let's denote by $\Ga^{st}$ \big (resp. $\wtilde{\Ga}^{st}$\big ), $s,t\in \{0,1\}$ the 
cuspidal divides obtained by performing $\Theta_s$,$\Theta_t$ splittings respectively at  
the $+$ double point and the vertical tangent point of $\Ga$
\big (resp.  of $\wtilde{\Ga}$\big ) in \fullref{fig-I-1}, without changing any other singular point of 
these divides. From lemma \ref{S3:L2}, we have splitting diagrams given in \fullref{fig-I-2}.
\\
Notice that $\Comp{\Ga ^{00}}=\Comp{\wtilde{\Ga}^{00}}$. In other words we have: 
%
\lem{\label{S4:SSI:L1} Let's denote 
$
\Delta^{0}=\left(\begin{array}{@{\, }c@{\, }}
				d^{0\bullet }\\
				d^{\bullet 0}
\end{array}\right )
$,
$
\Delta^{1}=\left(\begin{array}{@{\, }c@{\ }c@{\, }}
					d^{\bullet 1}	&	d^{1\bullet }	
					\end{array}\right )
$,
$
\wtilde \Delta^{0}=\left(\begin{array}{@{\, }c@{\, }}
				\wtilde d^{0\bullet }\\
				\wtilde d^{\bullet 0}
\end{array}\right )
$ and
$ 
\wtilde \Delta^{1}=\left(\begin{array}{@{\, }c@{\ }c@{\, }}
					\wtilde d^{\bullet 1}	&\wtilde d^{1\bullet }	
					\end{array}\right )
$. Then $\Comp{\Ga}=\cone\big (\Delta^{0},\Delta^{1}\big )$ and $\Comp{\wtilde{\Ga}}=\cone\big 
(\wtilde{\Delta}^{0},\wtilde{\Delta}^{1}\big )$:
\begin{equation}
	\begin{array}{c@{\ }c@{\ }c@{\ }c@{\ }c}
		\Comp{\Ga^{00}}\{3\}&		\displaystyle\mathop{\longrightarrow}^{\Delta^{0}}&
			\Comp{\Ga^{01}}\{4\}\oplus \Comp{\Ga^{10}}\{5\}&
			\displaystyle\mathop{\longrightarrow}^{\Delta^{1}}& 
			\Comp{\Ga^{11}}\{6\}\\
		\shortparallel\\
		\Comp{\wtilde \Ga^{00}}\{3\}&
		\displaystyle\mathop{\longrightarrow}^{\wtilde{\Delta}^{0}}&
		\Comp{\wtilde \Ga^{01}}\{4\}\oplus \Comp{\wtilde \Ga^{10}}\{5\}&
		\displaystyle\mathop{\longrightarrow}^{\wtilde{\Delta}^{1}}&
		\Comp{\wtilde\Ga^{11}}\{6\}
	\end{array}
	\label{DoubCompI}
\end{equation}
}	
%
Let's consider the ``creation and destruction'' morphisms (see remark \ref{S3:R2}):
$$
\begin{array}{c@{\, }c@{\, }c!{\quad }c@{\, }c@{\, }c!{\quad }c@{\, }c@{\, }c}
	\Comp{\Ga^{11}}\{6\}
	&\mathop{\longrightarrow }\limits^{\eta_1} 
	&\Comp{\Ga^{10}}\{5\} 
	&\Comp{\Ga^{01}}\{4\}
	&\mathop{\longrightarrow }\limits^{\wbar \eta_1}
	&\Comp{\wtilde \Ga^{10}}\{5\}
	&\Comp{\wtilde \Ga^{10}}\{5\}
	&\mathop{\longrightarrow }\limits^{\wtilde \varepsilon_1}
	&\Comp{\Ga^{01}}\{4\}\\
	\Comp{\wtilde \Ga^{11}}\{6\}
	&\mathop{\longrightarrow }\limits^{\wtilde \eta_1} 
	&\Comp{\wtilde \Ga^{10}}\{5\} 
	&\Comp{\wtilde \Ga^{01}}\{4\}
	&\mathop{\longrightarrow }\limits^{\wtilde {\wbar \eta}_1}
	&\Comp{\Ga^{10}}\{5\}
	&\Comp{\Ga^{10}}\{5\}
	&\mathop{\longrightarrow }\limits^{\varepsilon_1}
	&\Comp{\wtilde \Ga^{01}}\{4\}
\end{array}
$$
defined by  \fullref{fig-I-3}.
%
\lem{\label{S4:SSI:L2} The two sequences:
\begin{equation}
\begin {array}{c}
	0\longrightarrow \Comp{\Ga^{11}}\{6\}
		\mathop{\longrightarrow }\limits^{\eta_1}\Comp{\Ga^{10}}\{5\}
		\mathop{\longrightarrow }\limits^{\varepsilon_1}\Comp{\wtilde \Ga^{01}}\{4\}
		\longrightarrow 0\\
	0\longrightarrow \Comp{\wtilde \Ga^{11}}\{6\}
		\mathop{\longrightarrow }\limits^{\wtilde \eta_1}\Comp{\wtilde \Ga^{10}}\{5\}
		\mathop{\longrightarrow }\limits^{\wtilde \varepsilon_1}\Comp{\Ga^{01}}\{4\}
		\longrightarrow 0
\end {array}
\end{equation}
are exact and $d^{1\bullet}$, $\wtilde {\wbar \eta}_1$,$\tilde d^{1\bullet}$ and $\wbar \eta_1$ are respectively 
sections of $\eta_1$, $\varepsilon_1$, $\wtilde \eta_1$ and  $\wtilde \varepsilon _1$:
\begin{equation}
	\begin {array}{lll}
			d^{1\bullet}\eta_1=\id,&  \varepsilon_1\wtilde {\wbar \eta}_1=\id,& 
			\eta_1d^{1\bullet}+\wtilde {\wbar \eta}_1\varepsilon_1=\id+\eta_1d^{1\bullet}\wtilde {\wbar 
			\eta}_1\varepsilon_1\\
			\wtilde  d^{1\bullet}\wtilde \eta_1=\id,& \wtilde \varepsilon _1\wbar \eta_1=\id,& 
			\wtilde \eta_1\wtilde  d^{1\bullet}+ \wbar \eta_1\wtilde \varepsilon_1=\id+\wtilde \eta_1
			\wtilde  d^{1\bullet} \wbar \eta_1\wtilde \varepsilon_1
	\end {array}
\end{equation}
Moreover:
\begin{equation}
\eps_1 d^{\bullet 0}=\wtilde  d^{0\bullet }\qquad \text{and}\qquad 
\wtilde \eps_1\wtilde d^{\bullet 0}=d^{0\bullet }.
\label{rel}
\end{equation} 
}
%
\begin{proof} The morphisms $d^{1\bullet }$ and $\tilde d^{1\bullet }$ correspond to $T_7$ (see \ref{transfoT2}).
Then the result is an immediate consequence of remark \ref {S3:R2} (see also \fullref{fig3-1-1} and 
\fullref{fig3-1-3}).
\end{proof}
\begin{proof}[Proof of proposition \ref{S4:SSI:P1}.] Consider the diagram:
\begin{equation*}
\begin{array}{c@{\, }c@{\, }c@{\, }c@{\, }c}
	\Comp{\Ga^{00}}\{3\}&\mathop{\longrightarrow}\limits^{\Delta^0}& 
	\Comp{\Ga^{01}}\{4\}\oplus \Comp{\Ga^{10}}\{5\} &
	\mathop{\raise -5pt\hbox {$\stackrel{\dis \longrightarrow}{\curvearrowbotleft}$}}\limits^{\Delta^{1}}_{H}& 
	\Comp{\Ga^{11}}\{6\}
	\vspace{-5pt}
	\\
	\vspace{-5pt}
	\ \big \updownarrow \scriptstyle{\id}
	&
	&\scriptstyle{F}\big \downarrow \big \uparrow \scriptstyle{\wtilde F}
	&
	&\ \big \updownarrow \scriptstyle{0}\\
	\Comp{\wtilde \Ga^{00}}\{3\}&\mathop{\longrightarrow}\limits^{\wtilde \Delta^0}& 
	\Comp{\wtilde \Ga^{01}}\{4\}\oplus \Comp{\wtilde \Ga^{10}}\{5\} &
	\mathop{\raise -5pt\hbox {$\stackrel {\dis \longrightarrow}{\curvearrowbotleft}$}}\limits^{\wtilde 
	\Delta^{1}}_{\wtilde H}& 
	\Comp{\wtilde \Ga^{11}}\{6\}\\
\end{array}
\end{equation*}
\begin{gather*}
\tag*{where:} H=\left (\begin{array}{@{}c@{}}0 \\ \eta_1 \end{array}\right )\quad
\wtilde  H=\left (\begin{array}{@{}c@{}}0 \\ \wtilde \eta_1 \end{array}\right )\\
F=\left (\begin{array}{@{}cc@{}}
			0&\eps_1\\
			\wbar \eta_1+\wtilde \eta_1\wtilde d^{1\bullet }\wbar \eta_1 &\wtilde \eta_1\wtilde d^{\bullet 1}\eps_1 
\end{array}\right )\quad 
\wtilde F=\left (\begin{array}{@{}cc@{}}
			0&\wtilde  \eps_1\\
			\wtilde  {\wbar \eta}_1+\eta_1d^{1\bullet }\wtilde  {\wbar \eta}_1&\eta_1d^{\bullet 1}\wtilde  \eps_1
\end{array}\right ).
\end{gather*}
From lemma \ref{S4:SSI:L1} and \ref{S4:SSI:L2}, we have:
\begin{gather*}
	F\Delta^0=\wtilde  \Delta^0 \quad \wtilde \Delta^1 F=0\quad  \wtilde  FF=\id+H\Delta^1\quad \Delta^1H=\id\\
	\wtilde  F\wtilde  \Delta^0=\Delta^0\quad \Delta^1\wtilde  F=0\quad F\wtilde  F=\id+\wtilde  H\wtilde  
	\Delta^1\quad \wtilde  \Delta^1\wtilde  H=\id. 
\end{gather*}
Hense vertical arrows define a homotopy equivalence. From proposition \ref {S3:P4},  
$\Comp{\Ga}=\cone(\Delta^0,\Delta^1)$
and $\Comp{\wtilde \Ga}=\cone(\wtilde \Delta^0,\wtilde \Delta^1)$ have the same homology.
\end{proof}
\subsection{Invariance under type II move.}\label {S4:SS1} %
\begin{figure}[!ht]
  \centering
  \includegraphics{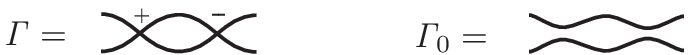}
	\caption{Type II move.}
	\label{fig-II-1}
\end{figure}

Let $\Ga$ and $\Ga_0$ be OMS--divides which differ only by a type II move (see \fullref {fig-II-1}).
%
\prop{\label{S4:SS1:P1} The complexes $\Comp{\Ga}$ and $\Comp{\Ga_0}$ have 
the same homology.}
Let $\Ga^{st}$, $s,t\in \{0,1\}$ be the cuspidal divides 
obtained by performing  $\Theta_s$ and $\Theta_t$ splittings respectively at the $+$ and the $-$ double points of 
$\Ga$ in \fullref {fig-II-1} without changing any other singular point. 
From lemma \ref{S3:L2} we have a splitting diagram given by \fullref {fig-II-2}.
We remark that $\Comp{\Ga^{00}}=\Comp{\Ga_0}$ and we have the following lemma:

\begin{figure}[!ht]
  \centering
  \includegraphics[width=\hsize]{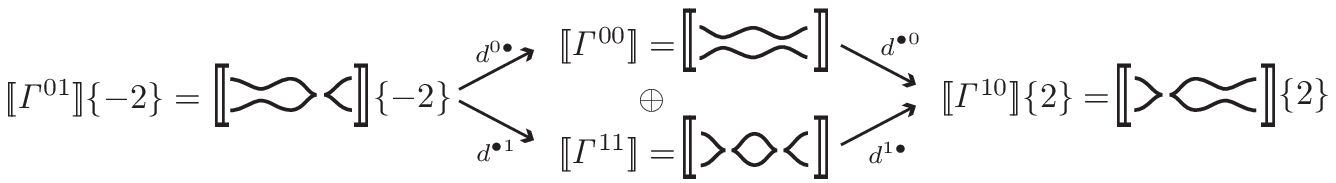}
	\caption{Splitting diagram for type II move.}
	\label{fig-II-2}
\end{figure}

\begin{figure}[!ht]
  \centering
  \includegraphics{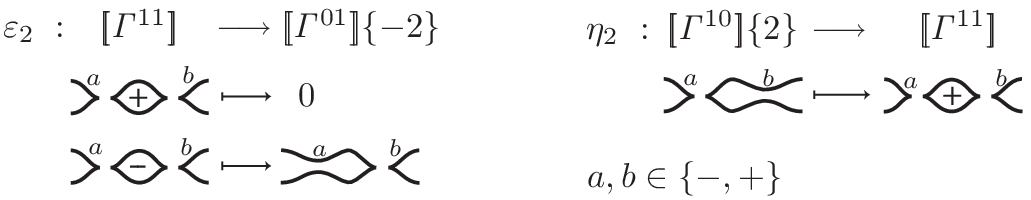}
	\caption {}
	\label{fig-II-3}
\end{figure}

%
\lem{\label{S4:SS1:L1} $\Comp{\Ga}=\cone(\Delta ^0,\Delta ^1)[-1]$ where:\
$$
\Comp{\Ga^{01}}\{-2\}
\stackrel{\Delta ^0}{\longrightarrow}
\Comp{\Ga^{00}}\oplus\Comp{\Ga^{11}}
\stackrel{\Delta ^1}{\longrightarrow}
\Comp{\Ga^{10}}\{2\}
\qquad 
\Delta ^0=\left(\begin{array}{@{}c@{}} d^{0\bullet } \\  d^{\bullet 1}		\end{array}	\right )
\quad 
\Delta ^1=\left (\begin{array}{@{}c@{\ }c@{}}  d^{\bullet 0}& d^{1\bullet } \end{array}\right ).
$$
}
%
Let's consider the ``destruction and creation'' morphisms of complexes (see remark \ref{S3:R2})
defined by \fullref{fig-II-3}.
%
\lem{\label{S4:SS1:L2} The sequence: \
$\dis 
0\longrightarrow \Comp{\Ga^{10}}\{2\}\mathop{\longrightarrow}^{\eta_2}
\Comp{\Ga^{11}}\mathop{\longrightarrow}^{\eps_2} \Comp{\Ga^{01}}\{-2\}\longrightarrow 0 
$\  
is exact and $d^{1\bullet }$ is a section of $\eta_2$, $d^{\bullet 1}$ 
a section of  $\eps_2$:
\begin {equation}
\eps_2 d^{\bullet 1}=\id,\qquad d^{1\bullet }\eta_2=\id,\qquad 
d^{\bullet 1}\, \eps_2+\eta_2\, d^{1\bullet }=\id + \eta_2\, d^{1\bullet }d^{\bullet 1}\eps_2.
\label{relationsII}
\end {equation}
}
%
\begin{proof} The morphism $d^{\bullet 1}$ corresponds to $T_2$ or $T_4$ and the morphism $d^{1\bullet }$ 
to $T_3$ or $T_5$ in (\ref{transfoT1}). The lemma is a direct consequence of remark \ref{S3:R2} (see also 
\fullref{fig3-1-1} and \fullref{fig3-1-2}).
\end{proof}
\begin{proof}[Proof of proposition \ref{S4:SS1:P1}.]
Consider the diagram:
$$ 
\begin{array}{c@{\ }c@{\ }c@{\ }c@{\ }c@{\ }c}
	&&\Comp{\Ga^{00}}=\Comp{\Ga_0}\\
	&&\Ajvert {15}{5}\scriptstyle{R}\left\uparrow\vrule height 10pt width 0pt\right \downarrow \scriptstyle{J}\\ 
	\Comp{\Ga^{01}}\{-2\}
	&\displaystyle\mathop{\raise -5pt\hbox {$\stackrel{\longrightarrow}{\curvearrowbotleft}$}}_{H^0}^{\Delta\! ^0}
	&\Comp{\Ga^{00}}\oplus \Comp{\Ga^{11}}
	&\displaystyle\mathop{\raise -5pt\hbox {$\stackrel{\longrightarrow}{\curvearrowbotleft}$}}_{H^1}^{\Delta\! ^1}
	&\Comp{\Ga^{10}}\{2\}\\	
\end{array}
\quad 
\begin{array}{cc}
J=\left(\begin{array}{@{}c@{}}\id \\ \eta_2 d^{\bullet 0}\end{array}\right )&
R=\left(\begin{array}{@{}cc@{}}\id & d^{0\bullet }\eps_2 \end{array}\right )\\
\\
H^0=\left(\begin{array}{@{}cc@{}}0 & \eps_2\end{array}\right )&
H^1=\left(\begin{array}{@{}c@{}}0 \\ \eta_2 \end{array}\right )
\end{array}
$$
From the previous two lemmas, $\Comp{\Ga_0}$ is a strong deformation retract of 
$\Comp{\Ga}=\cone(\Delta^0,\Delta^1)[-1]$ with retraction $R$, inclusion $J$, and homotopy $(H^0,H^1)$:
$$
R\Delta^0=0,\ \Delta^1J=0,\ H^0\Delta^0=\id,\ RJ=\id,\ JR=\id+\Delta^0H^0+H^1\Delta^1,
\ \Delta^1H^1=\id.
$$
Hence from proposition \ref {S3:P4} they have the same homology.
\end{proof}
\subsection{Invariance under type III move.}\label {S4:SS2}     %
In this section, we only consider the case of move III$_+$. The case of III$_-$ can be checked 
in a similar way: we have dual situations as is said in remark \ref{S3:R1} and in the proof of proposition 
\ref{S3:P2}. 

\begin{figure}[!ht]
  \centering
  \includegraphics{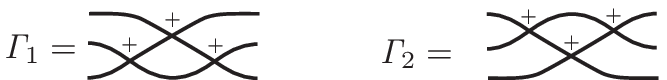}
	\caption{Type III$_+$ move.}
	\label{fig-III-1}
\end{figure}

Let $\Ga_1$ and $\Ga_2$ be OMS--divides which differ only by a 
type III$_+$ move (see \fullref{fig-III-1}).
%
\prop{\label{S4:SS2:P1}  The complexes $\Comp{\Ga_1}$ and $\Comp{\Ga_2}$ have 
the same homology.
}
%
Let $\Ga_1^{stu}$, $s,t,u\in \{0,1\}$ be the cuspidal divides obtained by performing
$\Theta_s$, $\Theta_t$ and $\Theta_u$ splittings at the double points shown on the figure of $\Ga_1$ 
(see \fullref{fig-III-1}). From lemma \ref{S3:L2}, we  have the following splitting diagram of $\Comp{\Ga_1}$ 
(see \fullref{fig-III-2}):
\begin{figure}[!ht]
  \centering
  \includegraphics[width=\hsize]{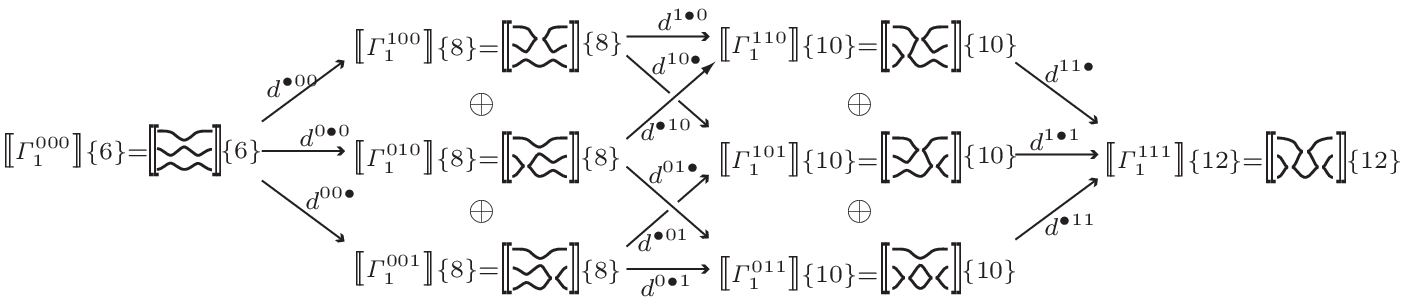}
	\caption{Splitting diagram for type III$_+$ move.}
	\label{fig-III-2}
\end{figure}
\\
Let's denote $\cC^0=\Comp {\Ga_1^{000}}\{6\}$,\ $\cC^3=\Comp {\Ga_1^{111}}\{ 12\}$,\
$\cC^{1}=\widehat \cC^{1}\oplus\widecheck \cC^{1}$ and $\cC^{2}=\widehat \cC^{2}\oplus	\widecheck \cC^{2}$\
\begin {gather*}
\tag*{where:}
\widehat \cC^{1}=\Comp {\Ga_1^{100}}\{8\}\oplus \Comp {\Ga_1^{010}}\{8\}\qquad 
\widecheck \cC^{1}=\Comp {\Ga_1^{001}}\{8\}\\
\widehat \cC^{2}=\Comp {\Ga_1^{110}}\{10\}\oplus \Comp {\Ga_1^{101}}\{10\}\qquad 
\widecheck \cC^{2}=\Comp {\Ga_1^{011}}\{10\}.
\end {gather*}
%
\lem{\label{S4:SS2:L1} $\Comp{\Ga_1}=\cone(\Delta^0,\Delta^1,\Delta^2)$\ where: \
$\dis 
\cC^0\mathop{\longrightarrow}^{\Delta^0}
\cC^1\mathop{\longrightarrow}^{\Delta^1}
\cC^2\mathop{\longrightarrow}^{\Delta^2}
\cC^3$ \ and:
\begin{equation*}
\Delta^0=\left (\begin{array}{@{}c@{}}
			d^{\bullet 00}\\
			d^{0\bullet 0}\\
			d^{00\bullet }
	\end{array}\right )
\qquad
\Delta^1=\left (\begin{array}{@{}ccc@{}}
			d^{1\bullet 0}&d^{\bullet 10}&0\\
			d^{10\bullet }&0&d^{\bullet 01}\\
			0&d^{01\bullet }&d^{0\bullet 1}
	\end{array}\right )
\qquad
\Delta^2=\left (\begin{array}{@{}ccc@{}} 
			d^{11\bullet }& d^{1\bullet 1}& d^{\bullet 11}
	\end{array}\right ).
\end{equation*}
}
%
Now we can use the same arguments as in \fullref{S4:SS1}.
The morphism $d^{0\bullet 1}$ (corresponding to $T_2$ and $T_4$) is injective 
and $d^{\bullet 11}$(corresponding to $T_3$ and $T_5$) is 
surjective.
Let $\eps_2: \widecheck \cC^{2}\to \widecheck \cC^{1}$ 
be the ``destruction'' morphism of complexes 
and $\eta_2: 	\cC^{3}\to \widecheck \cC^{2}$ 
be the ``creation'' morphism of complexes (see \fullref{fig-III-3}).
\begin{figure}[!ht]
	\centering
  \includegraphics{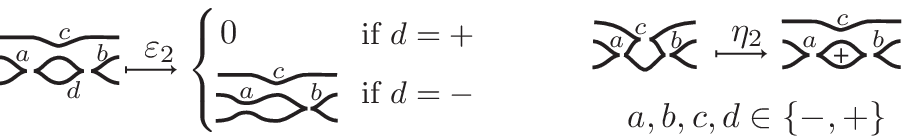}
	\caption{}
	\label{fig-III-3}
\end{figure}
%
\lem{\label{S4:SS2:L2} The sequence: $\dis 
0\longrightarrow \cC^{3}\mathop{\longrightarrow}^{\eta_2}\widecheck
\cC^{2}\mathop{\longrightarrow}^{\eps_2}\widecheck \cC^{1}\longrightarrow 0$
is exact and $d^{\bullet 11}$, $d^{0\bullet 1}$ are respectively sections of $\eta_2$ 
and $\eps_2 $:
\begin {equation}
d^{\bullet 11}\eta_2=\id,\quad \eps_2 d^{0\bullet 1}=\id,\quad
d^{0\bullet 1}\eps_2 +\eta_2 d^{\bullet 11}=\id +\eta_2 d^{\bullet 11}d^{0\bullet 1}\eps_2.
\label{relationsIII}
\end {equation}
}
%
\begin{proof} We have a similar situation as in lemma \ref {S4:SS1:L2}: $d^{0\bullet 1}$ corresponds to $T_2$ or 
$T_4$ and $d^{\bullet 11}$ corresponds to $T_3$ or $T_5$ 
in (\ref{transfoT1}). The lemma is a direct consequence of remark \ref{S3:R2} (see also \fullref{fig3-1-1}). 
\end{proof}
%
\lem{\label{S4:SS2:L3} Let $\delta=d^{\bullet 01}\eps_2 d^{01\bullet}:\Comp{\Ga_1^{010}}\{8\}\to 
\Comp{\Ga_1^{101}}\{10\}$.\\
The following sequence is a double-complex:
\begin{equation*}
(\cC^0,d^0)
\displaystyle\mathop{\longrightarrow}^{\widehat\Delta^0}
(\widehat \cC^{1},\widehat d^{1})
\displaystyle\mathop{\longrightarrow}^{\widehat\Delta ^1}
(\widehat \cC^{2},\widehat d^{2})\quad
\text{ where }\quad 
\widehat \Delta^0= \left (\begin{array}{@{}c@{}}
		 									d^{\bullet 00}\\ 
		 									d^{0\bullet 0}\\
											\end{array}\right )
											\quad 
\widehat \Delta^1= \left (\begin{array}{@{}cc@{}}
		 									d^{1\bullet 0}&d^{\bullet 10}\\
		 									d^{10\bullet }&\delta
											\end{array}\right ).
\end{equation*}
}
%
\proof Since $\Delta ^1\Delta^0=0$, we obtain from lemma \ref {S4:SS2:L2}:  
\begin{gather*}
\widehat \Delta^1\, \widehat \Delta^0
					=\left (\begin{array}{@{}c@{}}
		 									d^{1\bullet 0}d^{\bullet 00}+d^{\bullet 10}d^{0\bullet 0}\\
		 									d^{10\bullet }d^{\bullet 00}+\delta d^{0\bullet 0}
					\end{array}\right )=
					\left (\begin{array}{@{}c@{}}
		 									0\\
		 									d^{10\bullet }d^{\bullet 00}+d^{\bullet 01}\eps_2 d^{01\bullet}d^{0\bullet 0}
					\end{array}\right )\\
\tag*{and: }\dis d^{10\bullet }d^{\bullet 00}+d^{\bullet 01}\eps_2 d^{01\bullet}d^{0\bullet 0}=
d^{\bullet 01}d^{00\bullet }+d^{\bullet 01}\eps_2 d^{0\bullet 1}d^{00\bullet }=0.
\rlap{\hspace{1.6cm}\qedsymbol}
\end{gather*}
%
\lem{\label{S4:SS2:L4} Let $\widehat \cC=\cone( \widehat \Delta^0, \widehat \Delta^1)$. 
Then the complex $\widehat \cC$ is a strong deformation retract of $\Comp{\Ga_1}$ and so they have the 
same homology.
}
%
\begin{proof}
Consider the diagram:
$$
\begin{array}{ccccccc}
	\cC^{0} 
	&\hspace {-5pt}\stackrel {\widehat \Delta ^0}{\hbox to 35pt{\rightarrowfill }}\hspace {-20pt}
	&\widehat \cC^{1}
	&\hspace {-20pt}\stackrel {\widehat \Delta^1}{\hbox to 35pt{\rightarrowfill }}\hspace {-20pt}
	&\widehat \cC^{2}\\
	\left\updownarrow\vrule height 10pt width 0pt\right . \scriptstyle{\id}\hspace{-5pt}
	&&\scriptstyle{R^1}\left\uparrow\vrule height 10pt width 0pt\right \downarrow \scriptstyle{J^1}
	&&\scriptstyle{R^2}\left\uparrow\vrule height 10pt width 0pt\right \downarrow \scriptstyle{J^2}\\ 
	\cC^{0}
	&\displaystyle \mathop{\longrightarrow}^{\Delta^0}
	&\widehat \cC^{1}\oplus \widecheck \cC^{1}
	&\displaystyle\mathop{\raise -5pt\hbox {$\stackrel{\longrightarrow}{\curvearrowbotleft}$}}_{H^1}^{\Delta\! ^1}
	&\widehat \cC^{2}\oplus \widecheck \cC^{2}
	&\displaystyle\mathop{\raise -5pt\hbox {$\stackrel{\longrightarrow}{\curvearrowbotleft}$}}_{H^2}^{\Delta\! ^2}
	&\cC^{3}
\end{array}
$$ 
$$
\begin{array}{lll}
J^1=\left(\begin{array}{@{}cc@{}}\id &0\\ 0&\id\\ 
	0 & \eps_2 d^{01\bullet}\end{array}\right )
&R^1=\left(\begin{array}{@{}ccc@{}}\id & 0&0\\0&\id &0\end{array}\right )
&H^1=\left(\begin{array}{@{}ccc@{}}0 & 0&0\\0&0&0\\
0&0&\eps_2\end{array}\right )\\
J^2=\left(\begin{array}{@{}cc@{}}\id &0\\ 0&\id\\ 
\eta_2 d^{11\bullet}&\eta_2 d^{1\bullet
	1}\end{array}\right )
& R^2=\left(\begin{array}{@{}ccc@{}}\id & 0&0\\0&\id &d^{\bullet 01}\eps_2 \end{array}\right )
&H^2=\left(\begin{array}{@{}c@{}}0 \\0 \\ 
\eta_2 \end{array}\right )
\end{array}
$$
From lemma \ref {S4:SS2:L2}, we easily verify the relations:
\begin{gather*}
J^1\widehat \Delta ^0=\Delta ^0,\ \ J^2\widehat \Delta ^1=\Delta ^1J^1,\ \ \Delta ^2J^2=0,\ \
R^1\Delta ^0=\widehat \Delta ^0,\ \ R^2\Delta ^1=\widehat \Delta ^1R^1,\ \ R^1J^1=\id
\\ 
R^2J^2=\id,\ \ J^1R^1=\id +H^1\Delta ^1,\ \
J^2R^2=\id +H^2\Delta ^2+\Delta ^1H^1,\ 0=\id+\Delta^2H^2.
\end{gather*}
So downward arrows define an inclusion map, upward ones a retraction and $H_1,H_2$ homotopy maps of double 
complexes. We can apply proposition \ref {S3:P4}: $\widehat \cC$ is a deformation retract of $\Comp{\Ga_1}$
so they have the same homology. 
\end{proof}
%
\begin{proof}[Proof of proposition \ref{S4:SS2:P1}] 
Up to isotopy, the drawn part of $\Ga_1^{000}$ on \fullref{fig-III-2} is symmetric with respect to 
horizontal direction. Also the drawn part of $\Ga_1^{100}$ is symmetric to the one of $\Ga_1^{010}$ 
as well as the drawn part of $\Ga_1^{110}$ is symmetric to the one of $\Ga_1^{101}$. The morphisms 
$d^{\bullet 00}$ and $d^{0\bullet 0}$ (resp. $d^{10\bullet }$ and $d^{\bullet 10}$) clearly play symmetric roles.
Besides, the morphisms $\delta : \Comp{\Ga_1^{010}}\{8\}\to \Comp{\Ga_1^{101}}\{10\}$ and 
$d^{1\bullet 0}:\Comp{\Ga_1^{100}}\{8\}\to \Comp{\Ga_1^{110}}\{10\}$ also play symmetric roles.
Since the drawn parts of $\Ga_1$ and $\Ga_2$ in  \fullref{fig-III-1} also are symmetric with respect 
to horizontal direction, we deduce that twice the complexes $\Comp{\Ga_1}$ and 
$\Comp{\Ga_2}$ have the same homology as the complex $\widehat  \cC$.
\end{proof}
%
\subsection{Invariance under type IV moves.}\label {S4:SS4} %
\begin{figure}[!ht]
  \centering
  \includegraphics{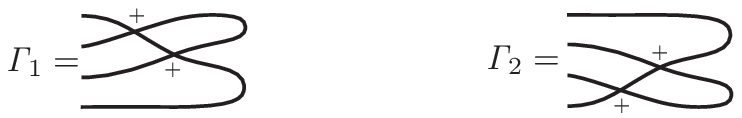}
	\caption{Type IV move.}
	\label{fig-IV-1}
\end{figure}

Let $\Ga_1$ and $\Ga_2$ be OMS--divides which differ only by a type $IV$ move 
(see \fullref{fig-IV-1}).
%
\prop{\label{S4:SS4:P1}  The complexes $\Comp{\Ga_1}$ and $\Comp{\Ga_2}$ have 
the same homology.
}
%
Before proving this proposition, we first introduce the following intermediate result.
%

\lem{\label{S4:SS4:L1} Let $\Ga$ and $\wtilde \Ga$ be two cuspidal divides 
which differ only in the following way:
\begin{figure}[!ht]
  \centering
  \includegraphics{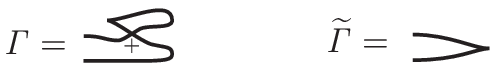}
	\caption{}
	\label{fig-IV-2}
\end{figure}
\\
(or symmetric situation with respect to horizontal direction).
Then $\Comp{\wtilde\Ga}\{6\}[1]$ is a strong deformation retract of $\Comp{\Ga}$.
}	
%
\begin{proof} With the same arguments as in the previous sections, using lemma \ref{S3:L2} we have a 
splitting diagram of $\Comp{\Ga}$ (\fullref{fig-IV-3}):
\begin{figure}[!ht]
  \centering
  \includegraphics[width=\hsize]{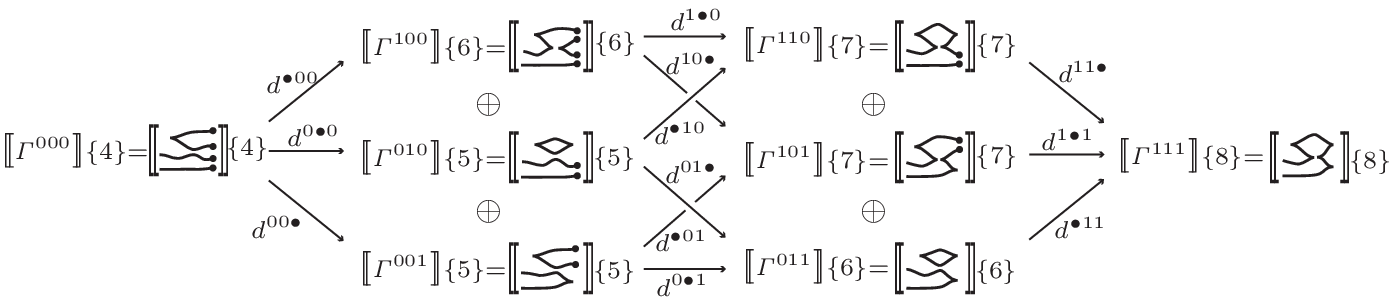}
	\caption{Splitting diagram of $\Comp{\Ga}$.}
	\label{fig-IV-3}
\end{figure}
\\
Let's denote:
\begin{gather*}
\cC^{0}=\Comp{\Ga^{000}}\{4\}\qquad 
\widecheck \cC^1=\Comp{\Ga^{001}}\{5\}
\qquad\quad
\widecheck \cC^2=\Comp{\Ga^{011}}\{6\}
\qquad
\cC^{3}=\Comp{\Ga^{111}}\{8\}
\\
\qquad
\widehat \cC^1=\Comp{\Ga^{100}}\{6\}\oplus \Comp{\Ga_1^{010}}\{5\}
\quad
\widehat \cC^2=\Comp{\Ga^{110}}\{7\}\oplus \Comp{\Ga_1^{101}}\{7\}.
\end{gather*}
Then $\Comp{\Ga}=\cone(\Delta^0,\Delta^1,\Delta^2)$ where: \
$\dis \cC^{0} \mathop{\longrightarrow}^{\Delta^0}\widehat \cC^1\oplus \widecheck \cC^1
\mathop{\longrightarrow}^{\Delta^1}\widehat \cC^2\oplus \widecheck \cC^2
 \mathop{\longrightarrow}^{\Delta^2}\cC^3
$
\begin{gather*}
\Delta^0=\left(\begin{array}{@{}c@{}}
	d^{\bullet 00}\\d^{0\bullet 0}\\ 
	d^{00\bullet }\end{array}\right )
=\left(\begin{array}{@{}c@{}}\Ajvert {17}{10} 
	\widehat \Delta^0\\ 
	d^{00\bullet }\end{array}\right )
\quad
\Delta^1=\left(\begin{array}{@{}ccc@{}}
		d^{1\bullet 0}&d^{\bullet 10}&0\\
		d^{10\bullet }&0&d^{\bullet 01}\\
		0&d^{01\bullet }&d^{0\bullet 1}\end{array}\right )
		=\left(\begin{array}{cc@{}}
		\widehat \Delta^1&U\\
		L&d^{0\bullet 1}\end{array}\right )
\\
\tag*{and }\Delta^2=\left(\begin{array}{@{}ccc@{}}d^{11\bullet }&d^{1\bullet 1}&d^{\bullet 11}\end{array}\right 
)=\left(\begin{array}{@{}cc@{}}\ \widehat \Delta^2 \ &d^{\bullet 11}\end{array}\right ).
\end{gather*}
Consider the creation / destruction morphisms of complexes (see \fullref{fig-IV-4}): 
$$
\begin{array}{@{}lll@{}}
\Comp {\Ga^{100}}\{6\}		\mathop{\to }\limits^{\tau }\		\Comp {\Ga^{000}}\{4\}
&\Comp {\Ga^{110}}\{7\}		\mathop{\to }\limits^{\eta_2}\		\Comp {\Ga^{010}}\{5\}
&\Comp {\Ga^{101}}\{8\}		\mathop{\to }\limits^{\eta_2}\		\Comp {\Ga^{011}}\{5\}
\\ 
\Comp {\Ga^{010}}\{5\}		\mathop{\to }\limits^{\sigma }\		\Comp {\Ga^{000}}\{4\}
&\Comp {\Ga^{101}}\{7\}		\mathop{\to }\limits^{\eta_1}\		\Comp {\Ga^{100}}\{5\}
&\quad \Comp {\wtilde \Ga}\{6\}	\mathop{\to }\limits^{\eta_1 }\		\Comp {\Ga^{001}}\{5\}
\\ 
&\Comp {\Ga^{101}}\{7\}		\mathop{\to }\limits^{\eta_2}\		\Comp {\Ga^{010}}\{5\}
&\Comp {\Ga^{001}}\{5\}		\mathop{\to }\limits^{\wbar \varepsilon_1}\	\Comp {\wtilde \Ga}\{6\}	
\end{array}
$$
\begin{figure}[!ht]
  \centering
  \includegraphics{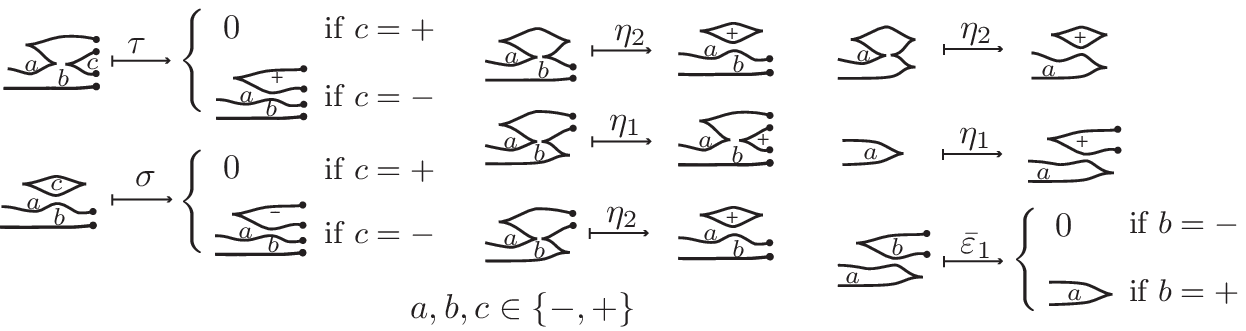}
	\caption{}
	\label{fig-IV-4}
\end{figure}
\\
Let's define \ $\sigma'=(\id + \tau d^{\bullet 00})\sigma : \Comp{\Ga^{010}}\{5\}\to \Comp{\Ga^{000}}\{4\}$ and: 
$$
\widehat H^1=\left(\begin{array}{@{}cc@{}}
		0&\eta_1\\
		\eta_2&\eta_2
\end{array}\right )
\qquad 
\widehat H^0=\left(\begin{array}{@{}cc@{}}
		\tau&\sigma'
\end{array}\right ).
$$ 
Then we have a short exact sequence \quad $\dis 
0\longrightarrow\widehat \cC^2
\mathop{\longrightarrow}^{\widehat H^1}
\widehat \cC^1
\mathop{\longrightarrow}^{\widehat H^0}
\cC^0\longrightarrow 0
$\quad such that $\widehat \Delta^1$ and $\widehat \Delta^0$ are sections of $\widehat H^1$
and $\widehat H^0$:
\begin{equation}
	\widehat \Delta^1\widehat H^1=\id,\quad \widehat H^0\widehat \Delta^0=\id,\quad 
	\widehat \Delta^0\widehat H^0+\widehat H^1\widehat \Delta^1=\id+
	\widehat H^1\widehat \Delta^1\widehat \Delta^0\widehat H^0.
	\label{RelationsIV}
\end{equation}
Moreover, for the compositions\ 
$\dis \Comp{\wtilde \Ga}\{6\}\mathop{\longrightarrow}^{\eta_1}
\widecheck \cC^1\mathop{\longrightarrow}^{d^{0\bullet 1}}\widecheck \cC^2$
\
and
\
$\dis \widecheck \cC^1\mathop{\longrightarrow}^{U}\widehat \cC^2 
\mathop{\longrightarrow}^{\widehat H^1}\widehat \cC^1
\mathop{\longrightarrow}^{L}\widecheck \cC^2
$
we have:
\begin{equation}
d^{0\bullet 1}\eta_1=0\qquad 
L\widehat H^1U=
d^{01\bullet}\eta_2d^{\bullet 01}=0.
\label{RelationsIV2}
\end{equation}
Let's define homotopies \ $H^0:\widehat \cC^{1}\oplus \widecheck \cC^{1}\to \cC^{0}$, \
$H^1: \widehat \cC^{2}\oplus \widecheck \cC^{2}\to \widehat \cC^{1}\oplus \widecheck \cC^{1}$\
and \ $H^2:\cC^{3}\to \widehat \cC^{2}\oplus \widecheck \cC^{2}$ \ by:
$$
H^0=\left(\begin{array}{@{}cc@{}}\widehat H^0 & 0\end{array}\right )
\qquad 
H^1=\left( \begin{array}{@{}cc@{}}
			\widehat H^1&
				\begin{array}{@{}c@{}} 0 \\ 
				 0\end{array}\\ 
			\begin{array}{@{}c@{\hskip 5pt}c@{}} 0 &  0\end{array}& 0
\end{array}\right )
\qquad H^2=
\left(\begin{array}{@{}c@{}}\!
		\begin{array}{c} 0 \\  0\end{array}\! \\ 
		 \eta_2\end{array}\right )
$$
and consider the injection $J$ and retraction $R$:
$$
\begin{array}{@{}c@{\hskip 5pt}c@{\hskip 5pt}c@{\hskip 5pt}c@{\hskip 5pt}c@{\hskip 5pt}c@{\hskip 5pt}c@{}}
	\\
	&&\ \ \Comp{\wtilde \Ga}\{6\}\\
	&&\Ajvert{12}{7}\scriptstyle{R}\left\uparrow\vrule height 10pt width 0pt\right \downarrow \scriptstyle{J}\\ 
	\cC^{0}
	&\displaystyle\mathop{\raise -5pt\hbox {$\stackrel{\longrightarrow}{\curvearrowbotleft}$}}_{H^0}^{\Delta\! ^0}
	&\widehat \cC^{1}\oplus \widecheck \cC^{1}
	&\displaystyle\mathop{\raise -5pt\hbox {$\stackrel{\longrightarrow}{\curvearrowbotleft}$}}_{H^1}^{\Delta\! ^1}
	&\widehat \cC^{2}\oplus \widecheck \cC^{2}
	&\displaystyle\mathop{\raise -5pt\hbox {$\stackrel{\longrightarrow}{\curvearrowbotleft}$}}_{H^2}^{\Delta\! ^2}
	&\cC^{3}
\end{array}
$$
$$
\begin{array}{c@{\ }c@{\ }c} 
	R & = & \left(\begin{array}{@{}ccc@{}} 
				 			0 &  0 & \wbar \eps_1 
					\end{array}\right )(\id+\Delta^0H^0)\\
	  & = &\Ajvert{15}{5}	\left(\begin{array}{@{}ccc@{}} 
						 \wbar \eps_1d^{00\bullet}\tau & 
						 \wbar \eps_1d^{00\bullet}\sigma' & 
						 \wbar \eps_1
					\end{array}\right )
\end{array}
\qquad
J=(\id+H^1\Delta^1)
			\left(\begin{array}{@{}c@{}}
					 0 \\ 
					 0 \\ 
					 \eta_1
			\end{array}\right )
		=\left (\begin{array}{@{}c@{}}
				  \eta_1 d^{\bullet 01}\eta_1 \\ 
				 	\eta_2 d^{\bullet 01}\eta_1 \\ 
					\eta_1
			\end{array}\right )
$$
(see \fullref{fig-IV-5} for retraction $R$).
\begin{figure}[!ht]
  \centering
  \includegraphics[width=\hsize]{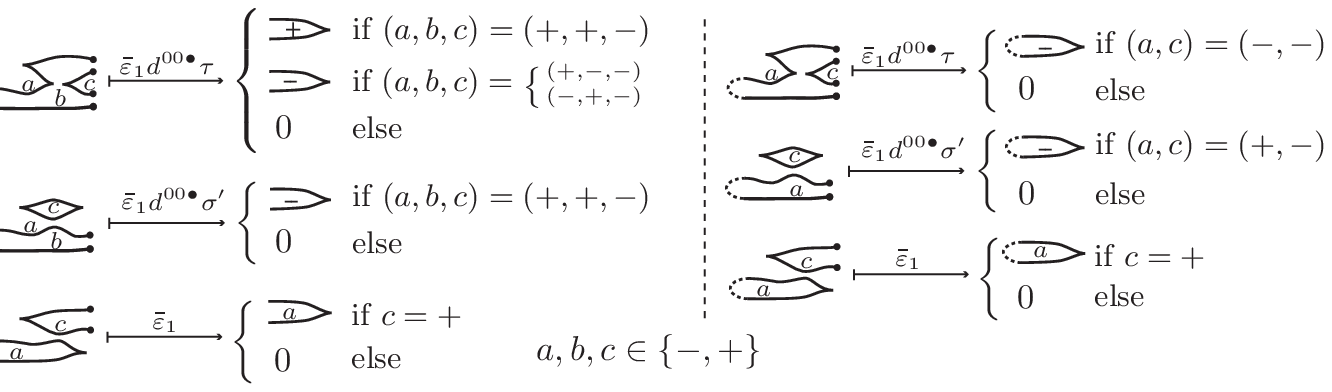}
	\caption{The retraction $R$.}
	\label{fig-IV-5}
\end{figure}
From (\ref {RelationsIV}) and (\ref{RelationsIV2}), we easily verify the relations:
\begin{gather*}
R \Delta ^0=0,\quad \Delta ^1J=0,\quad H^0J=0,\quad RH^1=0,\quad H^0H^1=0,\quad H^1H^2=0\\
RJ=\id,\ JR=\id +\Delta ^0H^0+H^1\Delta ^1,\ \id+\Delta^1H^1+H^2\Delta ^2=0,\ \id+\Delta^2H^2=0.
\end{gather*}
Then $\Comp{\wtilde \Ga}\{6\}[1]$ is a strong deformation retract of 
$\Comp{\Ga}=\cone(\Delta^0,\Delta^1,\Delta^2)$ from proposition \ref{S3:P4}.
\end{proof}
\begin{proof}[Proof of proposition \ref{S4:SS4:P1}]
Let $\Ga_1^s$ (resp. $\Ga_2^s$), $s\in \{0,1\}$ be the 
cuspidal divides obtained by performing $\Theta_s$ splittings at the left hand $+$ 
double point of $\Ga_1$ (resp. $\Ga_2$) in \fullref{fig-IV-1}. Let's also denote, according to lemma 
\ref{S4:SS4:L1}, the cuspidal divides $\wtilde \Ga^1_1$ (resp. $\wtilde \Ga^1_2$)
obtained by ``retracting'' $\Ga_1^1$ (resp. $\Ga^1_2$) (see \fullref{fig-IV-6}).
\begin{figure}[!ht]
	\centering
  \includegraphics{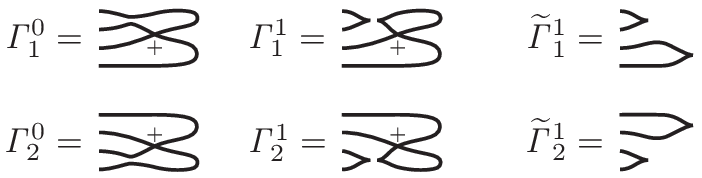}
	\caption{}
	\label{fig-IV-6}
\end{figure}
Notice that $\Comp{\Ga_1^0}=\Comp{\Ga_2^0}$ and $\Comp{\wtilde \Ga_1^1}=\Comp{\wtilde \Ga_2^1}$. 
From lemma \ref{S3:L2}, the differential $d_1$ on $\Comp{\Ga_1}$  
(resp. $d_2$ on $\Comp{\Ga_2$}) gives us the following cone:
$$
\Comp{\Ga_1}=\cone\big (\Comp{\Ga_1^0}\{2\}\mathop{\to}^{d_1^{\bullet}}\Comp{\Ga_1^1}\{4\}\big )
\quad 
\Big (\text{resp. }\ \Comp{\Ga_2}=\cone\big 
(\Comp{\Ga_2^0}\{2\}\mathop{\to}^{d_2^{\bullet}}\Comp{\Ga_2^1}\{4\}\big )\ \Big ).
$$
From lemma \ref{S4:SS4:L1}, there exist strong deformation retractions :
$$
\Comp{\Ga_1^1}\ \displaystyle\mathop{\raise -5pt\hbox {$\stackrel{\dis 
\longrightarrow}{\longleftarrow}$}}^{R_1}_{J_1}\ \Comp{\wtilde \Ga_1^1}\{6\}[1]
\qquad 
\Comp{\Ga_2^1}\ \displaystyle\mathop{\raise -5pt\hbox {$\stackrel{\dis
\longrightarrow}{\longleftarrow}$}}^{R_2}_{J_2}\ 
\Comp{\wtilde \Ga_2^1}\{6\}[1].
$$
So from proposition \ref{S3:P3}, $\cone (R_1d_1^{\bullet})$
(resp. $\cone (R_2d_2^{\bullet})$ ) is a strong deformation retract of $\Comp{\Ga_1}$ 
\big (resp. $\Comp{\Ga_2}$\big ). Hence it suffices to show that $R_1d_1^{\bullet}=R_2d_2^{\bullet}$. 
Let's consider the splitting diagram of $\Comp{\Ga_1^0}=\Comp{\Ga_2^0}$ as $\cone (D^0,D^1,D^2)$ 
(see \fullref{fig-IV-7}).
\begin{figure}[!ht]
  \centering
  \includegraphics[width=\hsize]{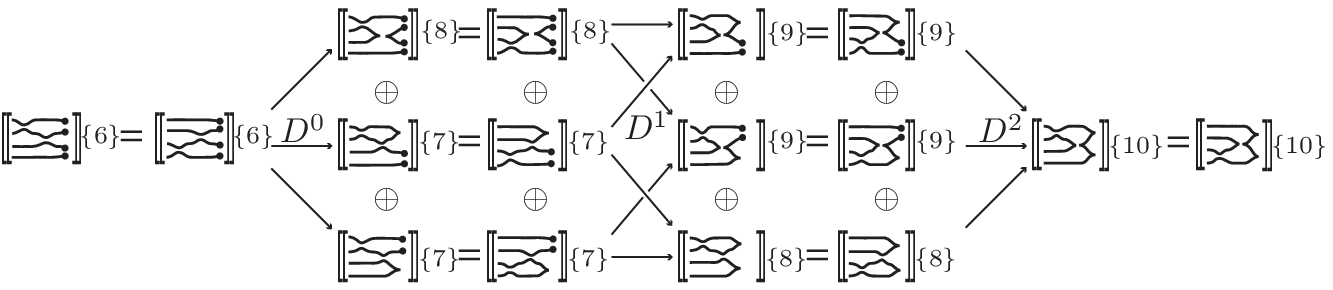}
	\caption{}
	\label{fig-IV-7}
\end{figure}
\\
Now $R_1d_1^{\bullet}$ and $R_2d_2^{\bullet}$ corresponds to the diagram of \fullref{fig-IV-8}.
\begin{figure}[!ht]
  \centering
  \includegraphics{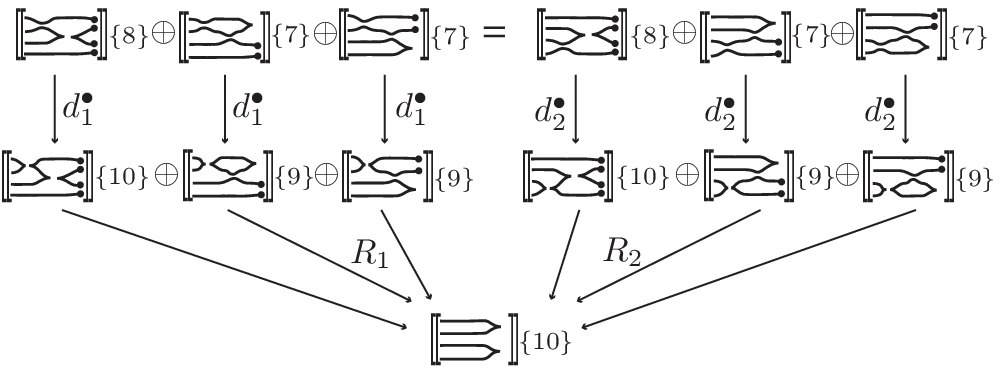}
		\caption{}
	\label{fig-IV-8}
\end{figure}
\\
By combining \fullref {fig-IV-5} with \fullref {fig3-1-1} we easily verify that 
$R_1d_1^{\bullet}=R_2d_2^{\bullet}$. 
\end{proof}
\subsection{Invariance under type V moves.}\label {S4:SS5} %
\begin{figure}[!ht]
  \centering
	\includegraphics{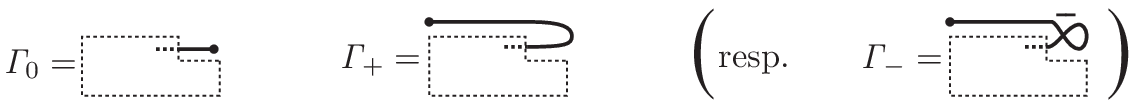}
	\caption{Type V moves.}
	\label{fig-V-1}
\end{figure}

Let $\Ga_0$ and $\Ga_+$ (resp. $\Ga_-$) be OMS--divides which differ only by a type $V_+$ (resp. type $V_-$) move 
(see \fullref{fig-V-1}).
%
\prop{\label{S4:SS5:P1}  The complexes $\Comp{\Ga_0}$, $\Comp{\Ga_+}$ and $\Comp{\Ga_-}$ have 
the same homology.
}
%
\begin{proof}
From lemma \ref{S3:L2}, we can see $\Comp{\Ga_+}$ as the cone of the surjective morphism:
\begin{figure}[!ht]
  \centering
  \includegraphics{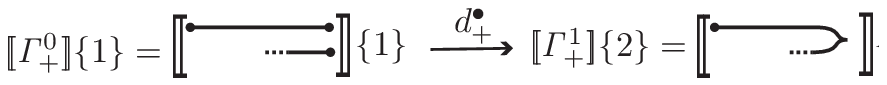}
	\caption{Splitting diagram for type V$_+$ move.}
	\label{fig-V-2}
\end{figure}
\\
Let's consider the creation / destruction morphisms
$\eta_1:\Comp{\Ga_+^1}\{2\}\to \Comp{\Ga_+^0}\{1\}$ (right inverse of $d_+^{\bullet}$),
$\wbar \eta_1:\Comp{\Ga}\to \Comp{\Ga_+^0}\{1\}$ and $\eps_1:\Comp{\Ga_+^0}\{1\}\to \Comp{\Ga}$
defined in \fullref{fig-V-3}.
\begin{figure}[!ht]
  \centering
  \includegraphics{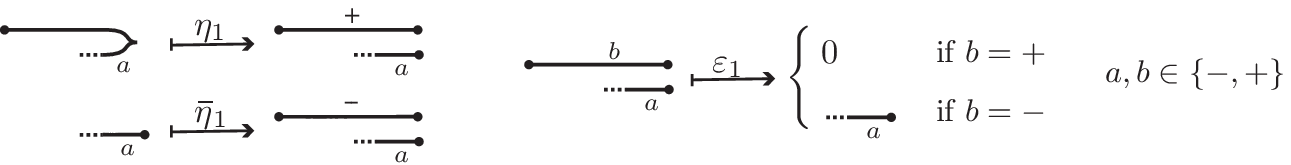}
	\caption{}
	\label{fig-V-3}
\end{figure}
Let $j=(\id+ \eta_1d^{\bullet})\wbar \eta_1$ and $r=\eps_1$. Then from the following diagram:
$$
\begin{array}{ccc}\Comp{\Ga_0}&\\
j\left  \downarrow\Ajvert{10}{0}\right \uparrow r\\
\Comp{\Ga_+^0}\{1\}&
\dis \mathop{\raise -5pt\hbox {$\stackrel{\longrightarrow}{\curvearrowbotleft}$}}_{\eta_1}^{d_+^{\bullet}}
&\Comp{\Ga_+^1}\{2\}
\end{array}
\qquad 
\begin{array}{c}
d_+^{\bullet}j=0,\quad rj=\id,\\
\\
jr=\id+\eta_1d_+^{\bullet},\quad d_+^{\bullet}\eta_1=\id
\end{array}
$$
we deduce using proposition \ref {S3:P4} that $\Comp{\Ga_0}$ is a strong deformation retract of 
$\Comp{\Ga_+}=\cone(d_+^{\bullet})$. So they have the same homology.
\\
On the other hand, from lemma \ref{S3:L2}, $\Comp{\Ga_-}=\cone(\Delta^0,\Delta^1)[-1]$:
$$\Comp{\Ga_-^{10}}\{-3\}\mathop{\to}^{\Delta^0}\Comp{\Ga_-^{00}}\{-1\}\oplus 
\Comp{\Ga_-^{11}}\{-2\}\mathop{\to}^{\Delta^1}\Comp{\Ga_-^{01}}\{3\}\ \
\Delta^0=\left (\begin{array}{@{}c@{}}d_-^{\bullet 0}\\
d_-^{1\bullet}\end{array}\right )\ \Delta^1=\left (\begin{array}{@{}c@{\ }c@{}}d_-^{0\bullet}&
d_-^{\bullet 1}\end{array}\right )
$$
\begin{figure}[!ht]
  \centering
  \includegraphics[width=\hsize]{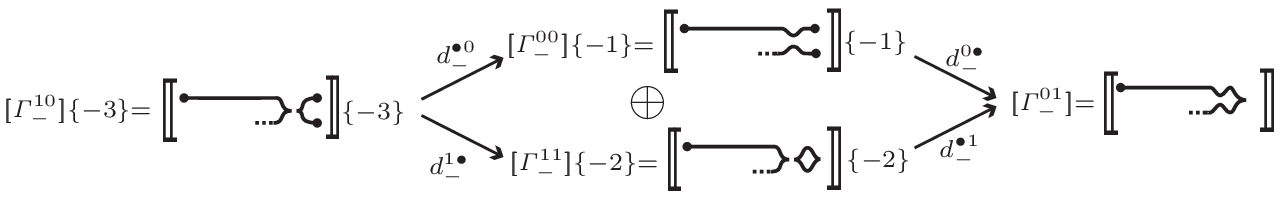}
\caption{Splitting diagram for type V$_-$ move.}
	\label{fig-V-4}
\end{figure}
\\
Let's consider the morphisms \ $\dis \Comp{\Ga_-^{00}}\{-1\}\mathop{\to}^{\tau}\Comp{\Ga_-^{10}}\{-3\}$, \
$\dis \Comp{\Ga_-^{11}}\{-2\}\mathop{\to}^{\sigma} \Comp{\Ga_-^{10}}\{-3\}$, \\
$\dis \Comp{\Ga_-^{01}}\mathop{\to}^{\eta_1}\Comp{\Ga_-^{00}}\{-1\}$\  and \ 
$\dis \Comp{\Ga_-^{01}}\mathop{\to}^{\eta_2}\Comp{\Ga_-^{11}}\{-2\}$ defined in \fullref{fig-V-5}.
\begin{figure}[!ht]
  \centering
  \includegraphics[width=\hsize]{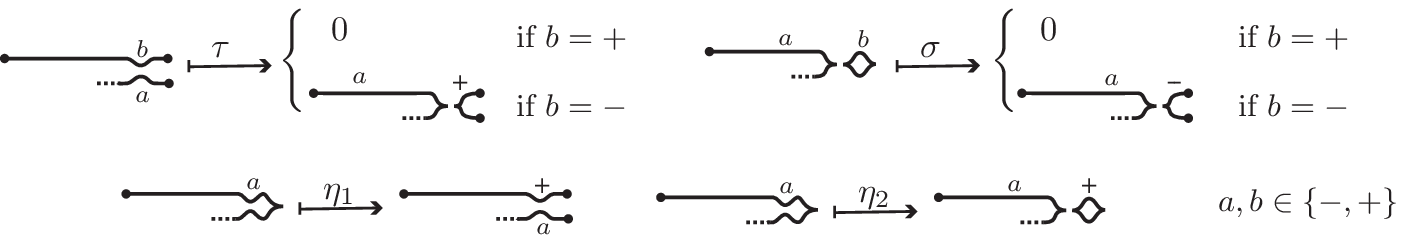}
	\caption{}
	\label{fig-V-5}
\end{figure}
\\
Then from the diagram:
$$
\begin{array}{ccccc}&&\Comp{\Ga_0}\\
&&{\scriptstyle J}\left  \downarrow\Ajvert{12}{0}\hskip -3pt\right \uparrow {\scriptstyle R}\\
\Comp{\Ga_-^{10}}\{-3\}&\dis \mathop{\raise -5pt\hbox 
{$\stackrel{\longrightarrow}{\curvearrowbotleft}$}}_{H^0}^{\Delta^0}&\Comp{\Ga_-^{00}}\{-1\}\oplus\Comp{\Ga_-^{11}}
\{-2\}&\dis \mathop{\raise -5pt\hbox 
{$\stackrel{\longrightarrow}{\curvearrowbotleft}$}}_{H^1}^{\Delta^1}&\Comp{\Ga_-^{01}}
\end{array}
$$
where $\dis J=\left (\begin{array}{@{}c@{}}\eta_1 \\ \eta_2\end{array}\right ) $, 
$\dis R=\left (\begin{array}{@{}cc@{}}0 & \wbar \eps_2\end{array}\right )$, 
$\dis H^0=\left (\begin{array}{@{}cc@{}}\tau & (\id +\tau d^{\bullet 0})\sigma \end{array}\right )$
and $\dis H^1=\left (\begin{array}{@{}c@{}}\eta_1\\ \eta_2\end{array}\right )$
we deduce that $\Comp{\Ga_0}$ is a strong deformation retract of $\Comp{\Ga_-}=\cone(\Delta^0,\Delta^1)[-1]$. 
They have the same homology.
\end{proof}
\subsection{Invariance under type VI moves.}\label {S4:SS6} %
\begin{figure}[!ht]
  \centering
  \includegraphics{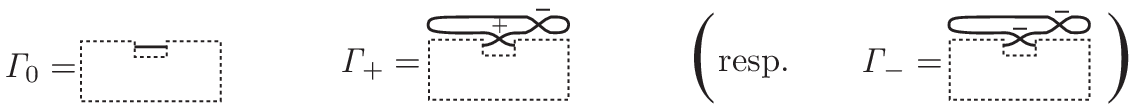}
	\caption{Type VI moves.}
	\label{fig-VI-1}
\end{figure}

Let $\Ga$ and $\Ga_+$ (resp. $\Ga_-$) be OMS--divides which differ only by a type $VI_+$ 
(resp. type $VI_-$) move (see \fullref{fig-VI-1}).
%
\prop{\label{S4:SS6:P1}  The complexes $\Comp{\Ga_0}$, $\Comp{\Ga_+}$ and $\Comp{\Ga_-}$ have 
the same homology.
}
%
We will break down the proof in two steps: the result is an immediate consequence of the following two lemmas.
%
\lem{\label{S4:SS6:L1}  Let $\Ga$ and $\wtilde \Ga$ be cuspidal divides defined by:
\begin{figure}[!ht]
  \centering
  \includegraphics{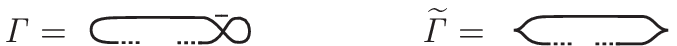}
	\caption{}
	\label{fig-VI-2}
\end{figure}
\\
Then $\Comp{\wtilde \Ga}$ is a strong deformation retract of $\Comp{\Ga}$.
}
%
\begin{proof} Let's apply lemma \ref{S3:L2} to $\Comp{\Ga}$. We have a splitting diagram :
\begin{figure}[!ht]
  \centering
  \includegraphics[width=\hsize]{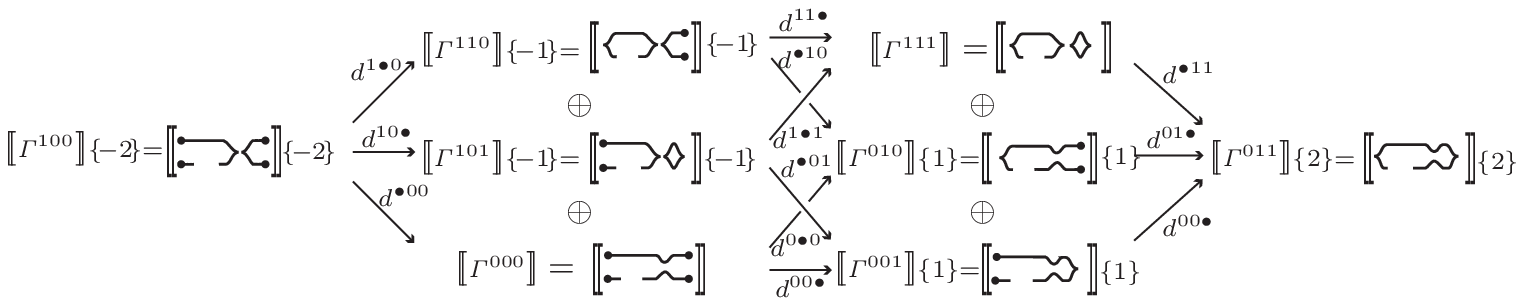}
	\caption{Splitting diagram of $\Comp{\Ga}$.}
	\label{fig-VI-3}
\end{figure}
\\
Let's denote :
\begin{gather*}
\cC^0=\Comp{\Ga^{100}}\{-2\}\qquad
\cC^1=\Comp{\Ga^{110}}\{-1\}\oplus \Big (\Comp{\Ga^{101}}\{-1\}\oplus \Comp{\Ga^{000}}\Big )
\\
\cC^2=\Comp{\Ga^{111}}\oplus\Big ( \Comp{\Ga^{010}}\{1\}\oplus\Comp{\Ga^{001}}\{1\}\Big )\qquad
\cC^3=\Comp{\Ga^{011}}\{2\}
\\
\Delta^0=\left(\begin{array}{@{}c@{}}
					d^{1\bullet 0}\\ 
					d^{10\bullet }\\ 
					d^{\bullet 00}
			\end{array}\right )
\qquad
\Delta^1=\left(\begin{array}{@{}ccc@{}}
					d^{11\bullet }&d^{1\bullet 1}&0\\
					d^{\bullet 10}&0&d^{0\bullet 0}\\
					0&d^{\bullet 01}&d^{00\bullet }
		\end{array}\right )
\qquad
\Delta^2=\left(\begin{array}{@{}ccc@{}}
					d^{\bullet 11}&d^{01\bullet }&d^{0\bullet 1}
			\end{array}\right ).
\end{gather*}
Then $\Comp{\Ga}=\cone(\Delta^0,\Delta^1,\Delta^2)[-1]$ where \ $\dis 
\cC^0\mathop{\longrightarrow}^{\Delta^0}\cC^1\mathop{\longrightarrow}^{\Delta^1}\cC^2
\mathop{\longrightarrow}^{\Delta^2}\cC^3$.
Let's consider the creation / destruction morphisms:
$$
\begin{array}{ccc}
	\dis \Ajvert{10}{7}\Comp{\Ga^{101}}\{-1\}\mathop{\longrightarrow}^{\sigma}  \Comp{\Ga^{100}}\{-2\} &
	\dis \Comp{\Ga^{111}}\mathop{\longrightarrow}^{\sigma}  \Comp{\Ga^{110}}\{-1\} &
	\Comp{\Ga^{011}}\{2\}\dis \mathop{\longrightarrow}^{\eta_2}\Comp{\Ga^{111}} 
	\\
	\dis \Ajvert{10}{7}\Comp{\Ga^{000}}\mathop{\longrightarrow}^{\tau}  \Comp{\Ga^{100}}\{-2\} &
	\Comp{\Ga^{010}}\{1\}\dis \mathop{\longrightarrow}^{\eta_1}\Comp{\Ga^{000}} &
	\Comp{\wtilde \Ga}\dis \mathop{\longrightarrow}^{\eta_1}\Comp{\Ga^{110}}\{-1\} 
	\\
	&
	\Comp{\Ga^{001}}\{1\}\dis \mathop{\longrightarrow}^{\eta_2}\Comp{\Ga^{101}}\{-1\} &  
	\Comp{\Ga^{110}}\{-1\}\dis \mathop{\longrightarrow}^{\wbar \eps_1} \Comp{\wtilde \Ga}
\end{array}		
$$
defined by:
\begin{figure}[!ht]
  \centering
  \includegraphics[width=\hsize]{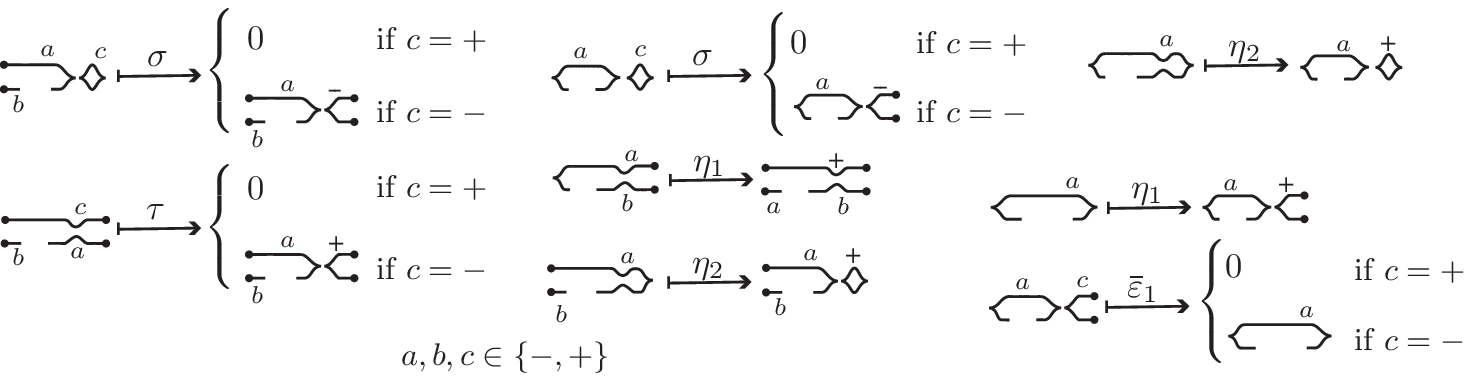}
	\caption{}
	\label{fig-VI-4}
\end{figure}
\\
We define homotopies $H^0:\cC^1\to \cC^0$,\  $H^1:\cC^2\to \cC^1$\ and\  $H^2:\cC^3\to \cC^2$:
$$
H^0=\left(\begin{array}{@{}ccc@{}}0 & (1+\tau d^{\bullet 00})\sigma & \tau\end{array}\right )
\
H^1=\left(\begin{array}{@{}ccc@{}}
		\sigma & 0 & 0\\
		\eta_2d^{00\bullet }\eta_1d^{\bullet 10}\sigma & \eta_2d^{00\bullet }\eta_1 & \eta_2\\
		\eta_1d^{\bullet 10}\sigma & \eta_1 & 0\end{array}\right )
\
H^2=\left(\begin{array}{@{}c@{}}\eta_2\\ 
		0\\0\end{array}\right ) 
$$ 
together with retraction and inclusion maps \ $R:\cC^1\to \Comp{\wtilde \Ga}$\  and \ 
$J:\Comp{\wtilde \Ga}\to \cC^1$:
$$
\begin{array}{r@{\ }c@{\ }l}
R&=&\left(\begin{array}{@{}ccc@{}} \wbar \eps_1 & 0 & 0\end{array}\right )(\id +\Delta^0H^0)\\
	\\
 &=&\left(\begin{array}{@{}ccc@{}} \wbar \eps_1 & \wbar \eps_1d^{1\bullet 0}\tau d^{\bullet 00}\sigma & 
		\wbar \eps_1d^{1\bullet 0}\tau \end{array}\right )
\end{array}
\ 
J=(\id +H^1\Delta^1)\left(\begin{array}{@{}c@{}} \eta_1\\ 
		0\\ 0 \end{array}\right )=
\left(\begin{array}{@{}c@{}} \eta_1\\ 
		\eta_2d^{00\bullet }\eta_1d^{\bullet 10}\eta_1\\ 
\eta_1d^{\bullet 10}\eta_1 \end{array}\right )
$$
such that we have the diagram:
$$
\begin{array}{cc@{}c@{}cccc}&&\Comp{\wtilde\Ga}\\
&&{\scriptstyle J}\left  \downarrow\Ajvert{12}{0}\hskip -3pt\right \uparrow {\scriptstyle R}\\
\cC^0&\dis \mathop{\raise -5pt\hbox 
{$\stackrel{\longrightarrow}{\curvearrowbotleft}$}}_{H^0}^{\Delta^0}&\cC^1&
\dis \mathop{\raise -5pt\hbox 
{$\stackrel{\longrightarrow}{\curvearrowbotleft}$}}_{H^1}^{\Delta^1}&\cC^2&
\dis \mathop{\raise -5pt\hbox 
{$\stackrel{\longrightarrow}{\curvearrowbotleft}$}}_{H^2}^{\Delta^2}&\cC^3
\end{array}
$$
We easily verify that :
\begin{gather*}
R \Delta ^0=0,\quad \Delta ^1J=0,\quad H^0J=0,\quad RH^1=0,\quad H^0H^1=0,\quad H^1H^2=0\\
RJ=\id,\ JR=\id +\Delta ^0H^0+H^1\Delta ^1,\ \id+\Delta^1H^1+H^2\Delta ^2=0,\ \id+\Delta^2H^2=0.
\end{gather*}
Then from proposition \ref{S3:P4}, $\Comp{\wtilde \Ga}$ is a strong deformation retract of 
$\Comp{\Ga}$.
\end{proof}
%
\lem{\label{S4:SS6:L2}  We have strong deformation retractions $r_+$ (resp. $r_-$) with 
injection $j_+$ (resp. $j_-$):
\begin{figure}[!ht]
  \centering
  \includegraphics[width=\hsize]{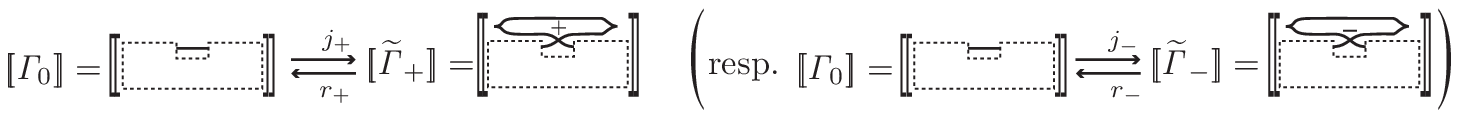}
	\caption{}
	\label{fig-VI-5}
\end{figure}
}
%
\begin{proof} From lemma \ref{S3:L2}, $\Comp{\wtilde \Ga_{+}}=\cone(d_+)$ and  $\Comp{\wtilde 
\Ga_{-}}=\cone(d_-)[-1]$ where:
\begin{figure}[!ht]
  \centering
  \includegraphics{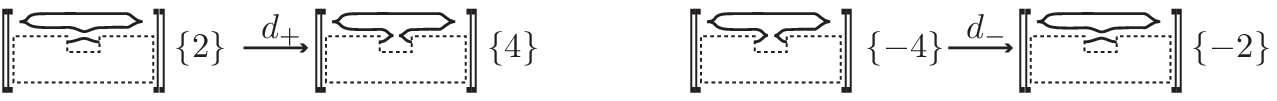}
	\caption{}
	\label{fig-VI-6}
\end{figure}
\\
For $\Comp{\Ga_{+}}$, using proposition \ref{S3:P4}, we deduce the strong deformation retraction 
from the following diagram of \fullref{fig-VI-7} since $d_+\eta_2=\id$,\ \ $r_+\eta_2=0$,\ \ $d_+j_+=0$,\ \ 
$r_+j_+=\id$ \ and \ $j_+r_+=\id+\eta_2d_+$.
\begin{figure}[!ht]
  \centering
  \includegraphics{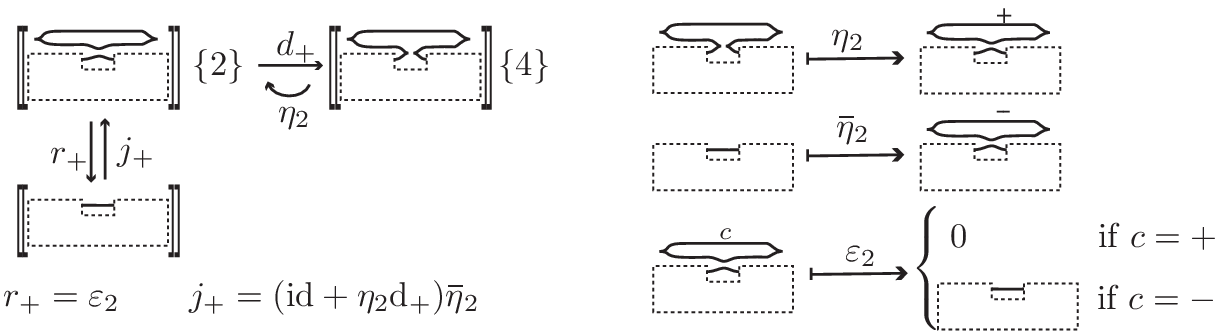}
	\caption{}
	\label{fig-VI-7}
\end{figure}
\\
Similarly for $\Comp{\Ga_{-}}$ from the diagram of \fullref{fig-VI-8}:
\begin{figure}[!ht]
  \centering
  \includegraphics{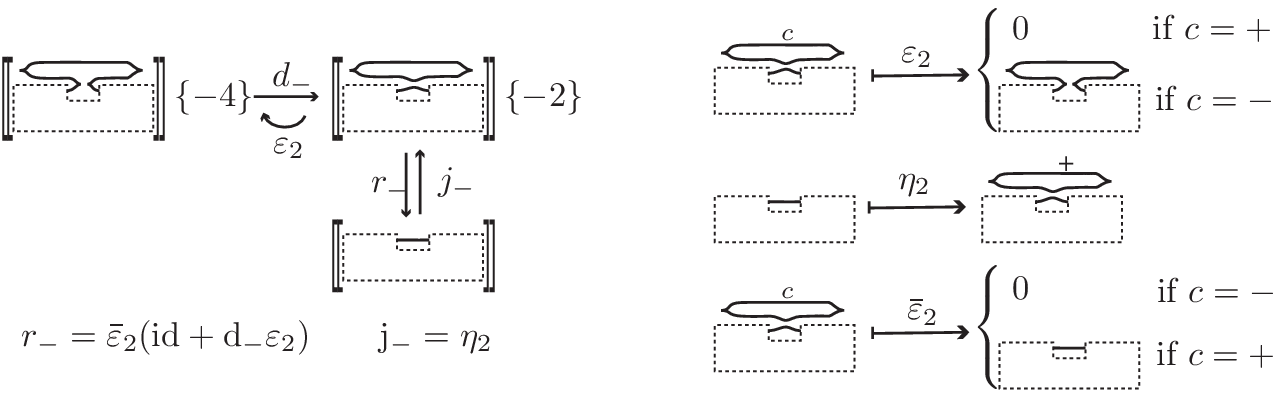}
	\caption{}
	\label{fig-VI-8}
\end{figure}
\\
since  $\eps_2d_-=\id$,\ $r_-d_-=0$,\ $\eps_2j_-=0$,\ $r_-j_-=\id$\ and \
$j_-r_-=\id+d_-\eps_2$.
\end{proof}

\end{document}